\documentclass[12pt,letterpaper,final,twoside,leqno]{amsart}
\usepackage{eucal, mathrsfs}
\usepackage{amsmath,amssymb,amsthm,
amscd,amsopn}
\usepackage{hyperref}
\usepackage{times}
\usepackage{xspace}
\usepackage{epsfig,epic,eepic,latexsym,color}
\usepackage[all]{xy}
\usepackage{paralist}
\usepackage{stmaryrd}
\usepackage{enumitem}
\usepackage{color}
\usepackage{soul}
\usepackage{mathtools}
\usepackage{rotating}

\catcode`~=11 \def\UrlSpecials{\do\~{\kern -.15em\lower .7ex\hbox{~}\kern .04em}} \catcode`~=13

\newcommand{\urlwofont}[1]{\urlstyle{same}\url{#1}}

\newcommand{\widepagestyle}{
\voffset=0in
\hoffset=0in
\marginparwidth=0.7in
\oddsidemargin=0in
\evensidemargin=0in
\textwidth=6.5in
\textheight=8.5in
\topmargin=0in
\headheight=0in
\headsep=0.2in
\footskip=0.5in
}

\newcounter{are-there-sections}
\makeatletter

\renewcommand\subsection{
  \renewcommand{\sfdefault}{pag}
  \@startsection{subsection}%
  {2}{0pt}{-\baselineskip}{.2\baselineskip}{\raggedright
    \sffamily\itshape\small\bfseries
  }}
\renewcommand\section{
  \renewcommand{\sfdefault}{phv}
  \@startsection{section} %
  {1}{0pt}{\baselineskip}{.2\baselineskip}{\centering
    \sffamily
    \scshape
    \bfseries
}}
\newcounter{lastyear}
\addtocounter{lastyear}{-1}

\newcommand\noin{\noindent}

\newtheoremstyle{bozont}{8pt}{10pt}%
     {\itshape}
     {}
     {\bfseries}
     {.}
     {.5em}
     {\thmname{#1}\thmnumber{ #2}\thmnote{ \rm #3}}
\newtheoremstyle{bozont-sf}{3pt}{3pt}%
     {\itshape}
     {}
     {\sffamily}
     {.}
     {.5em}
     {\thmname{#1}\thmnumber{ #2}\thmnote{ \rm #3}}
\newtheoremstyle{bozont-sc}{3pt}{3pt}%
     {\itshape}
     {}
     {\scshape}
     {.}
     {.5em}
     {\thmname{#1}\thmnumber{ #2}\thmnote{ \rm #3}}
\newtheoremstyle{bozont-remark}{8pt}{15pt}%
     {}
     {}
     {\scshape}
     {}
     {.5em}
     {\thmname{#1}\thmnumber{ #2.}\thmnote{\ \  ( #3)}}
\newtheoremstyle{bozont-def}{8pt}{15pt}%
     {}
     {}
     {\bfseries}
     {.}
     {.5em}
     {\thmname{#1}\thmnumber{ #2}\thmnote{ \rm #3}}
\newtheoremstyle{bozont-reverse}{3pt}{3pt}%
     {\itshape}
     {}
     {\bfseries}
     {.}
     {.5em}
     {\thmnumber{#2.}\thmname{ #1}\thmnote{ \rm #3}}
\newtheoremstyle{bozont-reverse-sc}{3pt}{3pt}%
     {\itshape}
     {}
     {\scshape}
     {.}
     {.5em}
     {\thmnumber{#2.}\thmname{ #1}\thmnote{ \rm #3}}
\newtheoremstyle{bozont-reverse-sf}{3pt}{3pt}%
     {\itshape}
     {}
     {\sffamily}
     {.}
     {.5em}
     {\thmnumber{#2.}\thmname{ #1}\thmnote{ \rm #3}}
\newtheoremstyle{bozont-remark-reverse}{3pt}{3pt}%
     {}
     {}
     {\sc}
     {.}
     {.5em}
     {\thmnumber{#2.}\thmname{ #1}\thmnote{ \rm #3}}
\newtheoremstyle{bozont-def-reverse}{3pt}{3pt}%
     {}
     {}
     {\bfseries}
     {.}
     {.5em}
     {\thmnumber{#2.}\thmname{ #1}\thmnote{ \rm #3}}
\newtheoremstyle{bozont-def-newnum-reverse}{3pt}{3pt}%
     {}
     {}
     {\bfseries}
     {}
     {.5em}
     {\thmnumber{#2.}\thmname{ #1}\thmnote{ \rm #3}}
\theoremstyle{bozont}    
\ifnum \value{are-there-sections}=0 {%
  \newtheorem{proclaim}{Theorem}
} 
\else {%
  \newtheorem{proclaim}{Theorem}[section]
} 
\fi
\newtheorem{thm}[proclaim]{Theorem}

\newtheorem{cor}[proclaim]{Corollary} 
 
\newtheorem{lem}[proclaim]{Lemma} 
\newtheorem{prop}[proclaim]{Proposition} 
\newtheorem{conj}[proclaim]{Conjecture}

\newtheorem{fact}[proclaim]{Fact}
\theoremstyle{bozont-sc}
\newtheorem{proclaim-special}[proclaim]{\specialthmname}

\theoremstyle{bozont-remark}

\newtheorem{rem}[proclaim]{Remark}

\newtheorem{example}[proclaim]{Example} 
 
\newtheorem{motivation}[proclaim]{Motivation} 
 
\newtheorem*{SubHeading*}{\SubHeadingName}%
\newtheorem{SubHeading}[proclaim]{\SubHeadingName}
\newtheorem{sSubHeading}[equation]{\sSubHeadingName}
\newenvironment{demo-r}[1]{\def\SubHeadingName{#1}\begin{SubHeading-r}}
  {\end{SubHeading-r}}%
\newenvironment{subdemo-r}[1]{\def\sSubHeadingName{#1}\begin{sSubHeading-r}}
  {\end{sSubHeading-r}} %
\newenvironment{demo*}[1]{\def\SubHeadingName{#1}\begin{SubHeading*}}
  {\end{SubHeading*}}%

\newtheorem{question}[proclaim]{Question}

\newtheorem{defn-thm}[proclaim]{Definition--Theorem}  

\theoremstyle{bozont-def}    
\newtheorem{defn}[proclaim]{Definition}
\newtheorem{notation}[proclaim]{Notation} 

\theoremstyle{bozont-reverse}    

\theoremstyle{bozont-reverse-sc}
\newtheorem{proclaimr-special}[proclaim]{\specialthmname}
{\def\specialthmname{#1}\begin{proclaimr-special}}%
{\end{proclaimr-special}}
\theoremstyle{bozont-remark-reverse}

\newtheorem{SubHeading-r}[proclaim]{\SubHeadingName}
\newtheorem{sSubHeading-r}[equation]{\sSubHeadingName}
\newtheorem{SubHeadingr}[proclaim]{\SubHeadingName}

\theoremstyle{bozont-def-newnum-reverse}    

\theoremstyle{bozont-def-reverse}

\newtheorem{newnumspecial}[proclaim]{\specialnewnumname}

\numberwithin{equation}{proclaim}

\numberwithin{figure}{section}

\newenvironment{enumerate-p}{
  \begin{enumerate}}
  {\setcounter{equation}{\value{enumi}}\end{enumerate}}
\newenvironment{enumerate-cont}{
  \begin{enumerate}
    {\setcounter{enumi}{\value{equation}}}}
  {\setcounter{equation}{\value{enumi}}
  \end{enumerate}}

\newlength{\swidth}
\makeatother

\DeclareMathAlphabet{\smallchanc}{OT1}{pzc}%
                                 {m}{it}
\DeclareFontFamily{OT1}{pzc}{}
\DeclareFontShape{OT1}{pzc}{m}{it}%
             {<-> s * [1.100] pzcmi7t}{}
\DeclareMathAlphabet{\mathchanc}{OT1}{pzc}%
                                 {m}{it}

\newcommand{\tf}{\tilde{f}}

\DeclareFontFamily{OMS}{rsfs}{\skewchar\font'60}
\DeclareFontShape{OMS}{rsfs}{m}{n}{<-5>rsfs5 <5-7>rsfs7 <7->rsfs10 }{}
\DeclareSymbolFont{rsfs}{OMS}{rsfs}{m}{n}
\DeclareSymbolFontAlphabet{\scr}{rsfs}

\newcommand{\sA}{\scr{A}}

\newcommand{\sE}{\scr{E}}
\newcommand{\sF}{\scr{F}}
\newcommand{\sG}{\scr{G}}
\newcommand{\sH}{\scr{H}}

\newcommand{\sL}{\scr{L}}
\newcommand{\sM}{\scr{M}}

\newcommand{\sO}{\scr{O}}

\newcommand{\sT}{\scr{T}}

\newcommand{\sX}{\scr{X}}

\newcommand{\bC}{\mathbb{C}}

\newcommand{\bH}{\mathbb{H}}

\newcommand{\bL}{\mathbb{L}}

\newcommand{\bP}{\mathbb{P}}
\newcommand{\bQ}{\mathbb{Q}}

\newcommand{\bZ}{\mathbb{Z}}

\newcommand{\cC}{\mathcal{C}}

\newcommand{\cU}{\mathcal{U}}

\newcommand{\fF}{\mathfrak{F}}

\newcommand{\fK}{\mathfrak{K}}

\newcommand{\fM}{\mathfrak{M}}

\newcommand{\fY}{\mathfrak{Y}}

\newcommand{\ofM}{\overline{\mathfrak{M}}}

\newcommand{\ofY}{\overline{\mathfrak{Y}}}

\DeclareMathOperator{\ks}{{ks}}

\DeclareMathOperator{\iks}{{iks}}

\DeclareMathOperator{\Aut}{Aut}

\DeclareMathOperator{\Hom}{Hom}

\DeclareMathOperator{\im}{{im}}

\DeclareMathOperator{\Isom}{Isom}

\DeclareMathOperator{\red}{red}

\DeclareMathOperator{\rk}{{rk}}

\DeclareMathOperator{\Spec}{{Spec}}

\DeclareMathOperator{\var}{{Var}}
\DeclareMathOperator{\Var}{{Var}}

\newcommand{\factor}[2]{\left. \raise 2pt\hbox{\ensuremath{#1}} \right/
        \hskip -2pt\raise -2pt\hbox{\ensuremath{#2}}}

\def\coh#1.#2.#3.{H^{#1}(#2,#3)}
\def\dimcoh#1.#2.#3.{h^{#1}(#2,#3)}
\def\hypcoh#1.#2.#3.{\mathbb H_{\vphantom{l}}^{#1}(#2,#3)}
\def\loccoh#1.#2.#3.#4.{H^{#1}_{#2}(#3,#4)}
\def\dimloccoh#1.#2.#3.#4.{h^{#1}_{#2}(#3,#4)}
\def\lochypcoh#1.#2.#3.#4.{\mathbb H^{#1}_{#2}(#3,#4)}
\def\ses#1.#2.#3.{0  \longrightarrow  #1   \longrightarrow 
 #2 \longrightarrow #3 \longrightarrow 0} 
\def\sesshort#1.#2.#3.{0
 \rightarrow #1 \rightarrow #2 \rightarrow #3 \rightarrow 0}
\def\dist#1.#2.#3.{  #1   \longrightarrow 
 #2 \longrightarrow #3 \stackrel{+1}{\longrightarrow} } 
\def\CDdist#1.#2.#3.{  #1   @>>>  #2  @>>>   #3 @>+1>> }  
\def\shortses#1.#2.#3.{0  \rightarrow  #1   \rightarrow 
 #2  \rightarrow   #3 \rightarrow  0}
\def\shortdist#1.#2.#3.{  #1   \rightarrow 
 #2  \rightarrow   #3 \stackrel{+1}{\rightarrow} }  
\def\ddist#1.#2.#3.#4.#5.#6.{\CD
#1 @>>> #2 @>>> #3 @>+1>> \\
@VVV @VVV @VVV \\
#4 @>>> #5 @>>> #6 @>+1>> 
\endCD}
\def\ddistun#1.#2.#3.#4.#5.#6.{\CD
#1 @>>> #2 @>>> #3 @>+1>> \\
@. @VVV @VVV  \\
#4 @>>> #5 @>>> #6 @>+1>> 
\endCD}
\def\Iff#1#2#3{
\hfil\hbox{\hsize =#1
\vtop{\noin #2}
\hskip.5cm 
\lower.5\baselineskip\hbox{$\Leftrightarrow$}\hskip.5cm
\vtop{\noin #3}}\hfil\medskip}
\newcommand{\union}\cup
\newcommand{\intersect}\cap
\newcommand{\Union}\bigcup
\newcommand{\Intersect}\bigcap
\def\myoplus#1.#2.{\underset #1 \to {\overset #2 \to \oplus}}

\newcommand{\of}{\overline{f}}

\newcommand{\oA}{\overline{A}}

\newcommand{\oX}{\overline{X}}
\newcommand{\oY}{\overline{Y}}

\widepagestyle

\title[Arakelov-Parshin rigidity of towers of curve fibrations]{Arakelov-Parshin rigidity of towers of curve fibrations}
\author{Zsolt Patakfalvi}

\address{Zsolt Patakfalvi, University of Washington, Department of Mathematics, Box 354350,
Seattle, WA 98195, U.S.A.
}
\email{pzs@math.washington.edu}
\urladdr{http://www.math.washington.edu/\~{}pzs}

\begin{document}

\begin{abstract}
Arakelov-Parshin rigidity is concerned with varieties mapping rigidly to the moduli stack  $\fM_h$ of canonically polarized manifolds. Affirmative answer for any class of maps implies  finiteness of the given class.

This article studies Arakelov-Parshin rigidity on an open subspace of $\fM_h$, on the locus $\fK\fF_h$ of iterated Kodaira fibrations. First, we prove rigidity for all complete curves mapping finitely onto $\fK\fF_h$.   Then, for generic affine curves mapping into $\fK\fF_h$, rigidity is shown when $\deg h =2$. The method used in the latter part is showing that the iterated Kodaira-Spencer map is injective. \end{abstract}

\maketitle

\tableofcontents


\section{Introduction}
\label{sec:introduction_rigiditiy}

According to Grothendieck's functor of points point of view, understanding all maps into a space is equivalent to understanding the space itself. This turns out to be a crucial observation for moduli spaces with no concrete description, e.g., for  the moduli stack $\fM_g$ of smooth curves  of genus $g$. Over the complex numbers, the first interesting class of maps into $\fM_g$ are finite morphisms from curves. One of the famous conjectures of Shafarevich, now a theorem of Parshin and Arakelov, concerns this class. 

\begin{notation}
\label{notation:shafarevich_classical} 
Fix an integer $g \geq 2$, a smooth (not necessarily projective) curve $U$, its smooth compactification $B$ and set $\Delta:= (B \setminus U)_{\red}$. 
\end{notation}

\begin{thm}[(Shafarevich Conjecture, \cite{Parshin_Algebraic_curves_over_function_fields}, \cite{Arakelov_Families_of_algebraic_curves})]
\label{thm:shafarevich_conjecture}
In the situation of Notation \ref{notation:shafarevich_classical},
\begin{enumerate}
\item \textbf{Finiteness (F):} there are finitely many finite maps $U \to \fM_g$,
\item \textbf{Hyperbolicity (H):}  if $U \to \fM_g$ is a finite map, then 
\begin{equation*}
 \deg \omega_B(\Delta) > 0 . 
\end{equation*}
\end{enumerate}
\end{thm}

Furthermore, $\textbf{\emph{(F)}}$  was decomposed into two parts by Arakelov and Parshin.

\begin{thm}[(Finiteness part of Shafarevich conjecture, \cite{Parshin_Algebraic_curves_over_function_fields}, \cite{Arakelov_Families_of_algebraic_curves})]
\label{thm:finitenes_shafarevich_conjecture}
In the situation of  Notation \ref{notation:shafarevich_classical},
\begin{enumerate}
\item \label{itm:finitenes_shafarevich_conjecture:boundedness} \textbf{Boundedness (B):} there are finitely many deformation equivalence classes of finite maps $U \to \fM_g$,
\item \label{itm:finitenes_shafarevich_conjecture:rigidity} \textbf{Rigidity (R):} every finite map $U \to \fM_g$ is rigid. That is, its deformation equivalence class contains only one element.
\end{enumerate}
\end{thm}

Note, that there exists also a number field version  of Theorem \ref{thm:shafarevich_conjecture} \cite{FG_EFA}. Faltings proved his famous theorem about finiteness of rational points of smooth curves of genus at least two using this number field version.

In the last two decades there has been an enormous progress in generalizing the, now classical, statements of Theorems \ref{thm:shafarevich_conjecture} and \ref{thm:finitenes_shafarevich_conjecture} to higher dimensions. In these generalizations, first   $\fM_g$ is replaced by its higher dimensional counterpart, the moduli space of canonically polarized manifolds $\fM_h$ with fixed Hilbert polynomial $h$ \cite{Viehweg_Quasi_projective_moduli}. Note that the compactification of the latter moduli space is an exciting ongoing project (c.f., \cite{Hacon_Kovacs_Classification_of_higher_dimensional_algebraic_varieties}, \cite{Kollar_Moduli_of_varieties_of_general_type}). In the most general form of higher dimensional Shafarevich conjecture, after replacing $\fM_g$ by $\fM_h$, usually arbitrary dimensional bases are allowed as well. The main subject of the present article is the generalization of $\textbf{\emph{(R)}}$. However, let us summarize first what we expect generally in higher dimensions.

\begin{notation}
\label{notation:shafarevich_higher_dimension}
Fix a numerical polynomial $h$, a manifold (i.e. smooth variety) $U$ and a smooth compactification $B$  of $U$ such that $\Delta:=(B \setminus U)_{\red}$ is a simple normal crossing divisor. 
\end{notation}

\begin{conj}[(Higher dimensional version of Shafarevich conjecture)]
\label{conj:shafarevich_conjecture_higher_dimensions}
In the situation of Notation \ref{notation:shafarevich_higher_dimension},
\begin{enumerate}
\item \textbf{(B):} there are finitely many deformation equivalence classes of maps $U \to \fM_h$,
\item \textbf{(R):} as it will be explained shortly, the obvious generalization fails, however there is \textbf{Viehweg's rigidity conjecture:} if $\nu : U \to \fM_h$ is a quasi-finite map such that for the induced family $f : X \to U$,  $\Omega_{X/U}$ is relatively ample, then $f$ is rigid, 
\item \textbf{(H) :} if $U \to \fM_h$ is quasi-finite, then $\omega_B(\Delta)$ is big (or equivalently $\kappa(B,\Delta)  =\dim B$).
\end{enumerate}
\end{conj}

There has been huge progress recently, especially regarding $\textbf{\emph{(B)}}$ and  $\textbf{\emph{(H)}}$. The former was entirely solved in \cite{Kovacs_Lieblich_Boundedness_of_families_of_canonically_polarized_manifolds}, while the latter has been shown up to $\dim B=3$ in \cite{Kebekus_Kovacs_The_structure_of_surfaces_and_threefolds_mapping_to_the_moduli_stack_of_canonically_polarized_varieties}. Without claiming any completeness, some of the major articles concerning $\textbf{\emph{(B)}}$ are \cite{Bedulev_Viehweg_Shafarevich_conjecture_for_surfaces_of_general_type_over_function_fields}, \cite{Viehweg_Zuo_On_the_isotriviality_of_families_of_projective_manifolds_over_curves}, \cite{Kovacs_Logarithmic_vanishing_theorems_and_Arakelov_Parshin_boundedness_for_singular_varieties} and \cite{Kovacs_Lieblich_Boundedness_of_families_of_canonically_polarized_manifolds}. A similar list for $\textbf{\emph{(H)}}$ is 
\cite{Miglorini_A_smooth_family_of_minimal_surfaces_of_general_type_over_a_curve_of_genus_at_most_one_is_trivial}, \cite{Kovacs_Smooth_families_over_rational_and_elliptic_curves}, \cite{Kovacs_On_the_minimal_number_of_singular_fibres_in_a_family_of_surfaces_of_general_type}, \cite{Kovacs_Algebraic_hyperbolicity_of_the_fine_moduli_spaces}, \cite{Viehweg_Zuo_Base_spaces_of_non_isotrivial_families_of_smooth_minimal_models}, \cite{Kovacs_Viehweg_s_conjecture_for_families_over_P_n}, \cite{Kovacs_Vanishing_theorems_boundedness_and_hyperbolicity_over_higher_dimensional_bases}, \cite{Kebekus_Kovacs_Families_of_canonically_polarized_varieties_over_surfaces}, \cite{Kebekus_Kovacs_The_structure_of_surfaces_and_threefolds_mapping_to_the_moduli_stack_of_canonically_polarized_varieties}, \cite{Kebekus_Kovacs_Families_of_varieties_of_general_type_over_compact_bases} and \cite{Jabbusch_Kebekus_Families_over_special_base_manifolds_and_a_conjecture_of_Campana}.

Similarly to $\textbf{\emph{(B)}}$ and $\textbf{\emph{(H)}}$, $\textbf{\emph{(R)}}$ has been drawing a  decent amount of attention. However, in contrast to the spectacular results in the other two parts, there is still very little known about $\textbf{\emph{(R)}}$. The basic reason is that for higher dimensional moduli spaces $U \to \fM_h$ being quasi-finite  is not enough assumption to obtain rigidity (see Example \ref{ex:not_rigid_max_var} or \cite[Theorem 0.3]{Viehweg_Zuo_Complex_multiplication}). Loosely speaking some stronger hyperbolicity or variational assumption (or both) is needed. Hence, the higher dimensional version of $\textbf{\emph{(R)}}$ is more of a question than a conjecture, and is as follows. For the precise definition of rigid morphisms and rigid families see Definition \ref{def:rigid}.

\begin{question}[$\textbf{\emph{(R)}}$ in higher dimensions]
\label{qtn:rigidity}
Which maps $\nu : U \to \fM_h$ are rigid? 
We are particularly interested in rigid maps such that for the associated family $f : X \to U$, all coverings $f' : X' \to U'$ of all quasi-finite pullbacks of $f$  still give rigid maps $U' \to \fM_{h'}$. 
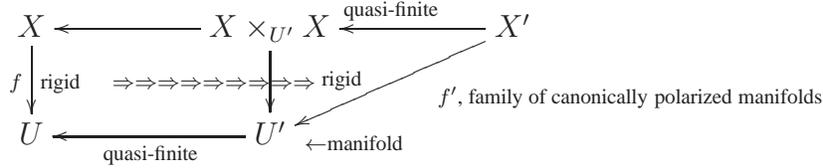
\begin{equation*}
\xymatrix{
X \ar[d]_f^{ \textrm{rigid} \quad \Rightarrow \Rightarrow \Rightarrow \Rightarrow \Rightarrow \Rightarrow \Rightarrow \Rightarrow \Rightarrow } & &  X \times_{U'} X \ar[ll] \ar[d]  & & X' \ar[ll]_{\textrm{quasi-finite}} \ar[dll]_(.65){\textrm{rigid}}^(.45){\qquad  \qquad  \qquad f' \textrm{, family of canonically polarized manifolds} } \\
U & & U' \ar[ll]^{\textrm{quasi-finite}} & \ar@{}[ul]^(0.10){\leftarrow \textrm{manifold}}
}
\end{equation*}
We call such maps \emph{stably rigid}. 
\end{question}

\begin{rem}
Stable rigidity seems to capture the philosophy of the rigidity condition of Theorem \ref{thm:finitenes_shafarevich_conjecture} (i.e. of the original Shafarevich conjecture). More precisely, stably rigid maps $U \to \fM_h$ from smooth curves are exactly finite maps by Proposition \ref{prop:stably_rigid}. We are expecting similar stability properties for any good rigidity condition in higher dimensions.
\end{rem}

Any positive answer to Question \ref{qtn:rigidity} yields $\textbf{\emph{(F)}}$ by $\textbf{\emph{(B)}}$. More precisely: 

\begin{thm} \cite[Theorem 1.6]{Kovacs_Lieblich_Boundedness_of_families_of_canonically_polarized_manifolds}
\label{thm:boundedness}
If for a fixed manifold $U$ and polynomial $h$, $\cC$ is a set of rigid maps $U \to \fM_h$, then $\cC$ is finite. 
\end{thm}

So far, there has been one answer in the literature to Question \ref{qtn:rigidity}, using the iterated Kodaira-Spencer map $\iks_{\nu}$ of Definition \ref{def:iterated_Kodaira_Spencer}. It is a notion motivated by Hodge Theory, and in case $\deg h =1$ it specializes to the ordinary Kodaira-Spencer map.

\begin{thm} \cite[Corollary 8.4]{Viehweg_Zuo_Discreteness_of_minimal_models_of_Kodaira_dimension_zero_and_subvarieties_of_moduli_stacks}, \cite[Theorem 4.14]{Kovacs_Strong_non_isotriviality_and_rigidity} 
\label{thm:iks_injective_rigid}
If $\nu : U \to \fM_h$ is a map from a smooth (not necessarily projective) curve such that $\iks_{\nu}$ is injective, then $\nu$ is rigid. 
\end{thm}

Although, Theorem \ref{thm:iks_injective_rigid} is a nice, general rigidity condition, it has certain unsatisfactory aspects. It is not known if it includes all stably rigid  families (c.f., \cite[Theorems 0.1 and 0.2]{Viehweg_Zuo_Complex_multiplication}). Also, we have no geometric understanding of the iterated Kodaira-Spencer map unless $\deg h=1$ or if $\nu$ corresponds to a family of hypersurfaces (c.f., \cite{Viehweg_Zuo_Complex_multiplication}). Being a notion motivated by Hodge Theory, understanding the iterated Kodaira-Spencer map is equivalent to understanding certain aspects of Torelli maps, which is very hard for higher dimensional canonically polarized manifolds. Therefore, we hope for easier ways to tackle Question \ref{qtn:rigidity} than via Theorem \ref{thm:iks_injective_rigid}.

\begin{rem}
Connections between variants of the Torelli problem and the Shafarevich conjecture have been around since the start of the subject. For example, $\textbf{\emph{(H)}}$ of Theorem \ref{thm:shafarevich_conjecture} over $\bC$ follows promptly from the Torelli theorem for curves and the boundedness of the weight one period domains. 
\end{rem}

\begin{rem}
\label{rem:gap}
There is a slight gap between \cite[Corollary 8.4]{Viehweg_Zuo_Discreteness_of_minimal_models_of_Kodaira_dimension_zero_and_subvarieties_of_moduli_stacks}, \cite[Theorem 4.14]{Kovacs_Strong_non_isotriviality_and_rigidity} and Theorem \ref{thm:iks_injective_rigid}. In the mentioned articles a different rigidity property is proven than what is used here (see Definition \ref{def:variationally_rigid}). However, the two rigidity properties turn out to be equivalent for one dimensional base by Proposition \ref{prop:variationally_rigid}. Note that during the course of  the proof of Proposition \ref{prop:variationally_rigid}, we also show higher dimensional Shafarevich conjecture for DM-curve bases in Lemma \ref{lem:stacky_shafarevich}.
\end{rem}

\subsection{Results of the paper}

As we have seen, the existing answers to Question \ref{qtn:rigidity} either concern the hypersurface case or have no geometric interpretations. In the present article we aim for geometric results in a case very different from that of hypersurfaces. We analyze the rigidity of towers of curve fibrations. More precisely, consider the following situation.

\begin{defn}
\label{def:tower}
A \emph{tower of curve fibrations} is a morphism $g: Y \to Z$ fitting in a commutative diagram
\begin{equation*}
\xymatrix{
Y=Y_n \ar@/^3pc/[rrrr]^g \ar[r]^{g_{n}} & Y_{n-1} \ar[r]^{g_{n-1}} & \dots  \ar[r]^{g_{2}} & Y_1 \ar[r]^{g_{1}} & Y_0= Z
},
\end{equation*}
where all schemes are varieties, $g_i$ are projective and the generic fibers of all $g_i$ are one dimensional and connected.
\end{defn}

\begin{notation}
\label{notation:tower}
Fix a tower of curve fibrations
\begin{equation}
\label{eq:tower_notation}
\xymatrix{
X=X_n \ar@/^3pc/[rrrr]^f \ar[r]^{f_{n}} & X_{n-1} \ar[r]^{f_{n-1}} & \dots  \ar[r]^{f_{2}} & X_1 \ar[r]^{f_{1}} & X_0= U
},
\end{equation}
such that $U$ is a smooth curve (not necessarily projective) over an algebraically closed field of characteristic zero, and $f_i$ are families of smooth curves with genera at least two. Let $\nu : U \to \fM_h$ be the moduli map associated to $f$. 
\end{notation}

\begin{rem}
Note that in the situation of Notation \ref{notation:tower}, if all $f_i$ are families of smooth curves of genus at least two, then $X_u$ is canonically polarized for all $u \in U$ by Proposition \ref{prop:rel_canonical_nef}. In particular, the map $\nu$ does exist.
\end{rem}

\begin{motivation}
Considering towers of curve fibrations is motivated partially by the following fact, which states that all families can be approximated in certain sense by towers of curve fibrations. Hence, we hope that in the long run, results about towers can be extended to general families.

By \cite[Corollary 5.10]{de_Jong_Families_of_curves_and_alterations} every family $ h: W \to U$ can be altered to a tower of curve fibrations as in Definition \ref{def:tower} such that $g_i$ are semi-stable families of curves. In other words, there is a commutative diagram 
\begin{equation*}
\xymatrix{
W \ar[d]^h  & & & & & Y  \ar[d]^g \ar[lllll]^{\textrm{generically finite, proper}} \\
U  & & & & & Z \ar[lllll]^{\textrm{generically finite, proper}},
}
\end{equation*}
with $g$ such a tower.
This fact would be even more promising if the answer to, the deliberately vaguely worded, Question \ref{qtn:alteration} was yes. It would mean, that every non-rigid family could be altered to a non-rigid tower of curve fibrations (using \cite[Corollary 1.4]{Hacon_McKernan_Boundedness_of_pluricanonical_maps_of_varieties_of_general_type}). Hence  stable rigidity could be determined by examining towers of curve fibrations only.
\end{motivation}

\begin{question}
\label{qtn:alteration}
If $g : W \to Z \times T$ is a deformation of the family $g_0 : W_{t_0} \to Z \times \{t_0\}$ of canonically polarized manifolds, is there then an alteration of the deformation $g$ into a  deformation of a tower of curve fibrations?
\end{question}

For our second motivation we need some definitions, which will be used extensively throughout the article.

\begin{defn}
If $f: X \to U $ is a family of canonically polarized manifolds, then the \emph{variation} of $f$ is $\var f := \dim ( \im \nu)$, where $\nu : U \to \fM_h$ is the associated moduli map. The family $f$ is called \emph{isotrivial} if $\var f = 0$.
\end{defn}

\begin{motivation}
In the situation of Notation \ref{notation:tower}, if $n=2$ and $\var f_i=i$, the question of rigidity is a special case of Viehweg's rigidity conjecture (e.g., \cite[Theorem 2]{Schneider_Complex_surfaces_with_negative_tangent_bundle}).
\end{motivation}

Now, we state the results of the paper. First for the case of a compact $U$, we have an almost full characterization of stable rigidity. Unfortunately, there is one possibility, mentioned in the last sentence of the theorem, which prevents us from giving a very concise answer. 

\begin{thm}
\label{thm:main_projective}
In the situation of Notation \ref{notation:tower}, if $U$ is projective, then 
\begin{enumerate}
 \item \label{itm:main_projective:rigid} if $\var f_i \geq 1$ for all $i$ then $\nu$ is rigid and
 \item otherwise there is a commutative diagram
 \begin{equation*}
  \xymatrix{
    X \ar[d]^f && X' \ar[ll]^{\textrm{\'etale}} \ar[d]^{f'} \cong W \times Y   \\
    U && U' \ar[ll]^{\textrm{\'etale, finite}}
  }
 \end{equation*}
 where $W \to U'$ is a family of canonically polarized manifolds, and $Y$ is a positive dimensional canonically polarized manifold (the map $W \times Y \to U'$ is the first projection composed with $W \to U'$). In particular, if $Y$ is not a rigid manifold, then $f'$ is not rigid. 
\end{enumerate}
\end{thm}



Theorem \ref{thm:main_projective} also implies rigidity of a class of compact curves on an open part of $\fM_h$. To make this more precise, we need to list a few more definitions.

\begin{defn}
\label{defn:iterated_Kodaira_fibration}
A \emph{family of iterated Kodaira fibrations} over a scheme $U$ is a commutative diagram as in \eqref{eq:tower_notation}, such that $f_i$ are families of smooth projective curves of genus at least two, and $\var \left( f_i|_{X_u} \right) \geq 1$ for each $i$ and $u \in U$. The Hilbert polynomial of such a family is $h(m) : = \chi(\omega_{X_u}^m)$, which is independent of the choice of $u \in U$. An \emph{iterated Kodaira fibration} is the special case of a  \emph{family of iterated Kodaira fibrations} with $U= \Spec k$.
The \emph{moduli space} $\widetilde{\fK \fF}_h$ \emph{of iterated Kodaira fibrations} is the stack, for which  
\begin{equation*}
 \left(\widetilde{\fK \fF}_h \right)_U=\left\{\parbox{185pt}{families of iterated Kodaira fibrations over $U$ with Hilbert polynomial $h$} \right\} , 
\end{equation*}
and morphisms of $\widetilde{\fK \fF}_h$ are the natural Cartesian diagrams with arrows at each level of the towers. This is a Deligne-Mumford stack of finite type over $k$ \cite[Proposition 2.9]{Patakfalvi_Fibered_stable_varieties}. Furthermore, there is a natural forgetful map $\pi : \widetilde{ \fK \fF }_h \to \fM_h$ forgetting all middle levels of the fibrations. Denote the image by $\fK \fF_h$. Then $\fK \fF_h$ is open in $\fM_h$, and the map $\pi : \widetilde{\fK \fF}_h \to \fK \fF_h$ is \'etale \cite[Theorem 1.2]{Patakfalvi_Fibered_stable_varieties}. 
\end{defn}

In fact,  by Theorem \ref{thm:proper}, $\fK \fF_h$ is open and closed in $\fM_h$. That is, it is a union of connected components. Also, it is easy to see that $\fK \fF_h$ is not empty for infinitely many values of $h$ by  Remark \ref{rem:iterated_Kodaia_fibrations_not_empty}. Furthermore, Theorem \ref{thm:main_projective} yields the following rigidity statement for these components.


\begin{cor}
\label{cor:main_projective}
If $\nu : U \to \fK \fF_h$ is a finite map from a smooth projective curve,  then $\nu$ is rigid.
\end{cor}

In the arbitrary base case, for  two level towers with maximal variations, the injectivity of the iterated Kodaira-Spencer map is shown as follows.

\begin{thm}
\label{thm:main_affine}
In the situation of Notation \ref{notation:tower}, if $n=2$ and $\var f_i = i$ (i.e. variations are maximal), then $\iks_{\nu}$ is injective (see Definition \ref{def:iterated_Kodaira_Spencer} for the definition of $\iks_{\nu}$). In particular, then $\nu$ is rigid by Theorem \ref{thm:iks_injective_rigid}.
\end{thm}


\begin{cor}
\label{cor:main_affine}
If $\nu : U \to \fK \fF_h$ is a finite map from a smooth curve, such that for a local lifting (i.e., $\tau$ is \'etale but not necessarily finite) 
\begin{equation}
\label{eq:main_affine:lifting}
\xymatrix{
U' \ar[r]^{\nu'} \ar[d]^{\tau} & \widetilde{\fK \fF}_h \ar[d] \\
U \ar[r]^{\nu} & \fM_h ,
}
\end{equation}
guaranteed by the \'etaleness of $\widetilde{\fK \fF}_h \to \fK \fF_h$, $\nu'$ corresponds to a tower as in Notation \ref{notation:tower} with $n=2$ and $\var f_i = i$, then $\nu$ is rigid. 
\end{cor}

An immediate corollary of Theorem \ref{thm:main_projective}, Theorem \ref{thm:main_affine} and Theorem \ref{thm:boundedness} is a finiteness statement.

\begin{cor}
Fixing a smooth curve $U$ and a polynomial $h$, there are finitely many
\begin{itemize}
\item maps $U \to \fK \fF_h$ with projective $U$,
\item maps $U \to \fM_h$ as in Notation \ref{notation:tower} with $\var f_i \geq 1$ and projective $U$  and
\item maps $U \to \fM_h$ as in Notation \ref{notation:tower} with $\var f_i =i$ and $n=2$.
\end{itemize}
\end{cor}

\begin{rem}
In Notation \ref{notation:tower}, singularities are not allowed. Still, in the statements of Theorem \ref{thm:main_affine} and Corollary \ref{cor:main_affine}, singularities are virtually allowed over finitely many points of $U$. Indeed, since the assumptions of these statements are not sensible to the restriction of the base, finitely many points of $U$ can be disregarded. However, even in these two statements, there cannot be singularities over generic $u \in U$. Note that similar sensibility to singularities can be observed in  Theorem \ref{thm:finitenes_shafarevich_conjecture}.\ref{itm:finitenes_shafarevich_conjecture:rigidity}, which fails if maps $U \to \ofM_g$ are allowed. So, we suspect that, for example, in Theorem \ref{thm:main_affine} singular fibers of $f_2$ could have not been allowed over generic $u \in U$. 
\end{rem}

%

\subsection{Organization of the paper}

The following two sections are concerned with preparations. In section \ref{sec:basic_concepts}, the basic definitions and statements left out from Section \ref{sec:introduction_rigiditiy}, to avoid technicalities there,  are collected. Section \ref{sec:positivity} is a short account on the results used in the article about the positivity of the relative canonical sheaves. Then in Section \ref{sec:compact}, Theorem \ref{thm:main_projective} and Corollary \ref{cor:main_projective} is proven. Along doing so, some facts about the moduli theory of families of canonically polarized manifolds is collected. In Section \ref{sec:relative_properness}, we prove the properness of $\widetilde{\fK \fF}_h \to \fM_h$, needed to deduce Corollary \ref{cor:main_projective} from theorem \ref{thm:main_projective}. Section \ref{sec:arbitrary base} is entirely devoted to the proof of Theorem \ref{thm:main_affine}. 
Finally, in Section \ref{sec:variationally_rigid} we fill the small gap mentioned in Remark \ref{rem:gap}.

\subsection{Notation} 

We work over an algebraically closed field $k$ of characteristic zero. 
All schemes are of finite type and separated over $k$ unless otherwise stated.  For a smooth projective curve $C$, $g(C)$ denotes its genus. A manifold is a smooth variety. A variety is an integral, separated scheme of finite type over $k$. A curve is a variety of dimension one. A global normal crossing divisor is defined Zariski locally by $\prod f_i^{n_i}$ where $f_i$ are regular elements and $n_i$ are positive integers. A canonically polarized manifold is a projective manifold $Z$ with ample $\omega_Z$. The Hilbert polynomial of a canonically polarized manifold $Z$ is $h(n):=\chi(\omega_Z^n)$. The Kodaira and log Kodaira dimensions of a variety $Z$ or a pair $(Z,\Delta)$ are denoted by $\kappa(Z)$ and $\kappa(Z, \Delta)$, respectively. For a line bundle $\sL$, its Iitaka-Kodaira dimension is denoted by $\kappa(\sL)$. We say, the 
variation of a family $g : Y \to Z$ is maximal if $\var g = \dim Z$. A vector bundle $\sE$ on $Y$ is ample over an open set $U$, if there is an ample line bundle $\sL$ and a homomorphism $\sL^{\oplus N} \to S^n(\sE)$ to some symmetric power of $\sE$, which is a surjection over $U$. $\sE$ is ample if it is ample over $Y$. We denote by $\fM_g$ and $\fM_h$ the moduli stacks of smooth projective curves of genus g and canonically polarized manifolds of Hilbert-polynomial $h$, respectively. Note that $\deg h$ is the dimension of the varieties parametrized by $\fM_h$.

\subsection{Acknowledgments}

I would like to thank my advisor, S\'andor Kov\'acs, for the fantastic guidance throughout my PhD studies. I would also like to thank J\'anos Koll\'ar for the useful remarks concerning a preprint version of the paper and Max Lieblich for discussing Section \ref{sec:variationally_rigid} with me.

\section{Basic concepts}
\label{sec:basic_concepts}

Here we collect some basic definitions and constructions mentioned in Section \ref{sec:introduction_rigiditiy}, which being slightly technical were omitted from there. We start with showing that $\fK \fF_h$ is not empty for infinitely many values of $h$.

\begin{rem}
\label{rem:iterated_Kodaia_fibrations_not_empty}
To see that $\fK \fF_h$ is not empty for infinitely many values of $h$, note that there are non-isotrivial smooth curve fibrations $p :Z  \to Y$ over projective curves such that both $Y$ and the fibers of $p$ have genus at least two \cite[V.14]{Barth_Peters_Van_de_Ven_Compact_complex_surfaces}. Furthermore, for infinitely many values of the fiber genus can be achieved. Define then 
\begin{equation*}
X_i:= 
\left\{
\begin{matrix}
\underbrace{Z \times_Y Z \times_Y \dots \times_Y Z }_{\textrm{$i-1$ times}} & \textrm{if } i \geq 2 \\
Y & \textrm{if } i=1 \\
\Spec k & \textrm{if } i=0 .  \\
\end{matrix}
\right.
\end{equation*}
and $f_i : X_i \to X_{i-1}$ the natural projection maps onto the first $i-1$ factors (or the adequate structure maps if $i=1$ and $2$). This yields a commutative diagram as \eqref{eq:tower_notation} which is an iterated Kodaira fibration by the choice of $p$. For different choices of $p$ we obtain infinitely many values of $h$, because the leading coefficient of $h$ is $\frac{K_{X_n}^n}{n!}$, and furthermore $K_{X_n}^n \geq \prod_{i=1}^n f_{n,i}^* K_{X_i/X_{i-1}} = g^n $, where $f_{n,i} : X_n \to X_i$ are the natural maps and $g$ is the genus of the fibers of $p$. 
\end{rem}

Next we give the precise definition of rigidity.

\begin{defn}
\label{def:rigid}
A morphism $f : X \to Y$ between Deligne-Mumford stacks is \emph{rigid}, if for every deformation $f' : X \times S \to Y$ over a smooth irreducible curve $S$, $f_s'=f$ for every $s \in S$. A family $f : X \to U$ of canonically polarized manifolds is \emph{rigid}, if the associated moduli map $U \to \fM_h$ is rigid.
\end{defn}

Next we show the promised example about why maximal variation does not imply rigidity for higher dimensional fibers.

\begin{example}
\label{ex:not_rigid_max_var}
Consider two non-isotrivial  families $f : S \to C$ and $g :T \to D$ of smooth projective curves of genera at least two with smooth projective bases  of genera at least two. Such families exist (e.g., \cite[Section V.14]{Barth_Peters_Van_de_Ven_Compact_complex_surfaces}). Consider $f \times g : S \times T \to C \times D$. It is a family of canonically polarized surfaces over $C \times D$. Moreover, since by \cite[Corollary 1.4]{Hacon_McKernan_Boundedness_of_pluricanonical_maps_of_varieties_of_general_type}, from a fixed variety there are only finitely many dominant maps onto varieties of general type up to birational equivalence, the restriction of $S \times T \to C \times D$ to $\{c \} \times D$ or $C \times \{d\}$ are non-isotrivial for any $c \in C$ and $d \in D$. So, fix any $d \in D$. Then $S \times T \to C \times D$ is a non-trivial deformation of the non-isotrivial family $S \times T_d \to C \times \{d\} \cong C$ of canonically polarized manifolds. However, since $\dim C=1$, here non-isotrivial means 
having maximal variation.  So, maximal variation does not imply rigidity in case of higher dimensional 
fibers.
\end{example}

The next proposition was promised after the statement of Question \ref{qtn:rigidity} and justifies the introduction of  stable rigid maps.

\begin{prop}
\label{prop:stably_rigid}
A map $\nu : U \to \fM_g$, for $g \geq 2$, from a smooth curve is stably rigid, if and only if it is finite.
\end{prop}

\begin{proof}
From Theorem \ref{thm:finitenes_shafarevich_conjecture}, using that non-isotriviality is stable under pulling back and taking cover (\cite[Corollary 1.4]{Hacon_McKernan_Boundedness_of_pluricanonical_maps_of_varieties_of_general_type}),  follows the backwards direction. To see the forward direction, assume $\nu$ is not finite, i.e., the associated family $f : X \to U$ is isotrivial. Then by Lemma \ref{lem:etale_trivialization} (with setting $Y:=X$, $U:=U$, $S:= \Spec k$), there is a finite \'etale cover $T \to U$, such that there is a diagram
\begin{equation*}
\xymatrix{
X \times_U T \ar[d] \ar@{<->}[r]^{\cong} & T \times F \ar[d] \\
T \ar@{=}[r] & T 
}
\end{equation*}
for some smooth projective curve $F$ of genus at least two. Then by deforming $F$, we get a deformation of $X \times_U T \to T$. That is, $f$ is not stably rigid.
\end{proof}

The rest of the section is devoted to the definition of the iterated Kodaira Spencer map. It is the main object of Section \ref{sec:arbitrary base}. 

\begin{defn}
\label{def:iterated_Kodaira_Spencer}
If $g : Y \to Z$ is a proper, smooth morphism of relative dimension $n$ over a smooth base, then for $1 \leq p \leq n$, by \cite[Exercise II.5.16]{Hartshorne_Algebraic_geometry} $\wedge^p \sT_Y$ has a filtration 
\begin{equation*}
0=\sF^p_0 \subseteq \sF^p_1 \subseteq \dots \subseteq \sF^p_{p} \subseteq \sF^p_{p+1}=\wedge^p \sT_Y 
\end{equation*}
 by locally free sheaves such that the induced quotients are 
\begin{equation}
\label{eq:iterated_Kodaira_Spencer:quotient}
\factor{\sF^p_{i+1}}{\sF^p_{i}} \cong (g^* \wedge^i \sT_Z) \otimes  (\wedge^{p-i} \sT_{Y/Z}) .
\end{equation}
Consider then the short exact sequences 
\begin{equation}
\label{eq:iterated_Kodaira_Spencer:short_exact}
\xymatrix{
0 \ar[r] 
& \wedge^p \sT_{Y/Z}  \ar[r] 
& \sF^p_2 \ar[r]   
& g^* \sT_Z \otimes \wedge^{p-1} \sT_{Y/Z}  \ar[r] & 0
}.
\end{equation}
Tensor these with $g^* \sT_Z^{\otimes n-p}$ to get the exact sequences
\begin{equation}
\label{eq:iterated_Kodaira_Spencer}
\xymatrix{
0 \ar[r] 
& g^* \sT_Z^{\otimes n-p} \otimes \wedge^p \sT_{Y/Z}  \ar[r] 
& g^* \sT_Z^{\otimes n-p} \otimes \sF^p_2     \\
&  \ar[r] 
& g^* \sT_Z^{\otimes n-p+1} \otimes \wedge^{p-1} \sT_{Y/Z}  \ar[r] & 0.
}
\end{equation}
Denote by $\rho_p$ the edge maps
\begin{equation*}
\rho_p : \sT_Z^{\otimes (n-p+1)} \otimes R^{p-1} g_*( \wedge^{p-1} \sT_{Y/Z}) \to \sT_Z^{\otimes (n-p)} \otimes R^{p} g_*  (\wedge^p \sT_{Y/Z})
\end{equation*}
obtained by applying higher pushforwards to  \eqref{eq:iterated_Kodaira_Spencer}. Then the Kodaira-Spencer map 
\begin{equation*}
\ks_g : \sT_Z \to R^1 g_* \sT_{Y/Z}
\end{equation*}
of $g$ is the edge map of \eqref{eq:iterated_Kodaira_Spencer:short_exact} for $p=1$. The  iterated Kodaira-Spencer map 
\begin{equation*}
\iks_g : \sT_Z^{\otimes n} \to R^n g_* (\wedge^n \sT_{Y/Z})
\end{equation*}
of $g$ to be $\rho_{n} \circ \dots \circ \rho_1$. We also define the $i$-th iterated Kodaira-Spencer map 
\begin{equation*}
\iks_g^i : \sT_Z^{\otimes n} \to \sT_Z^{\otimes (n-i)} \otimes  R^i g_* (\wedge^i \sT_{Y/Z})
\end{equation*}
by $\rho_{i} \circ \dots \circ \rho_1$. In particular, if $g$ is a family of canonically polarized manifolds and $\nu : Z \to \fM_h$ is the associated moduli map, then $\iks_{\nu} : = \iks_g$.
\end{defn}

\begin{rem}
In the case when $\dim Z=1$, $\sF^p_2 = \wedge^p \sT_Y$.
\end{rem}

\begin{rem}
There is another way to define $\iks_g$. It is the composition of the $n$ times product of $\ks_g$ and of the wedge product:
\begin{equation*}
\xymatrix{
\sT_Z^{\otimes n } \ar@/^2pc/[rrr]^{\iks_g} \ar[rr]_(0.4){\ks_f \otimes \dots \otimes \ks_f} & & (R^1 g_* \sT_{Y/Z})^{\otimes n } \ar[r]_{\wedge} & R^n g_* ( \wedge^n \sT_{Y/Z})
}
\end{equation*}
The equivalence of the two definitions can be proven using \v{C}ech or Dolbeault cohomology.
\end{rem}

\section{Positivity properties of the relative canonical sheaf}
\label{sec:positivity}

In this section certain positivity results are collected, some of which have already been used, and others will be used frequently later on. First, a statement about the relative canonical sheaves of a family of canonically polarized manifolds. 

\begin{prop}
\label{prop:rel_canonical_nef}
If $f : X \to B$ is a family of canonically polarized manifolds with $B$ smooth, projective, then $\omega_{X/B}$ is nef.
\end{prop}

\begin{proof}
It is known that $f_* \omega_{X/B}$  is a nef vector bundle (e.g., \cite[Theorem 4.1]{Viehweg_Weak_positivity}). Then, since $\omega_{X/B}$ is relatively ample, there is some $n>0$ such that $\omega_{X/B}^{\otimes n}$ is relatively globally generated. That is, there is a surjection $f^* f_* (\omega_{X/B}^{\otimes n}) \to \omega_{X/B}^{\otimes n}$, which shows the nefness of $\omega_{X/B}$.
\end{proof}

Next, another statement about the pushforwards of tensor powers of the relative canonical sheaf.

\begin{lem} \cite[Proposition 3.4]{Viehweg_Zuo_Base_spaces_of_non_isotrivial_families_of_smooth_minimal_models} 
\label{lem:positivity_pushforwards}
If $f : X \to B$ is a family of canonically polarized manifolds with $B$ smooth, projective and $\var f = \dim B$, then for any $\nu >1$, for which $f_* ( \omega_{X/B}^{\nu} ) \neq 0$, $f_*( \omega_{X/B}^{\nu})$ is ample with respect to the open subset $U \subseteq B$, where the moduli map $B \to \fM_h$ is quasi-finite.
\end{lem}

\begin{cor}
\label{cor:rel_canonical_sheaf_big}
In the situation of Lemma \ref{lem:positivity_pushforwards}, $\omega_{X/B}$ is ample with respect to $f^{-1} U$.
\end{cor}

\begin{proof}
Since $\omega_{X/B}$ is relatively ample, $\omega_{X/B}^n$ is relatively globally generated for $n \gg 0$. Choose such an $n$. Then there is a surjection 
\begin{equation*}
\omega_{X/B} \otimes f^* f_*( \omega_{X/B}^n) \to \omega_{X/B} \otimes \omega_{X/B}^n \cong \omega_{X/B}^{n+1},
\end{equation*}
which yields the statement of the lemma using Proposition \ref{prop:rel_canonical_nef} and that relatively ample nef line bundle tensored with the pullback of an ample vector bundle over $U$ is ample over $f^{-1}U$.
\end{proof}

\begin{cor}
\label{cor:rel_canonical_Kodaira_dim}
If $f : X \to B$ is a family of canonically polarized manifolds of dimension $n$ with $B$ smooth and projective, then $\kappa(\omega_{X/B})= \var f + n$.
\end{cor}

\begin{proof}
Let $\nu : B \to \fM_h$ be the moduli map. One can construct a commutative diagram
\begin{equation*}
\xymatrix{
& \tilde{X} \ar[dl]^{\xi} \ar[rrr]^{\zeta} \ar[dd]^{\tilde{f}} & & & Y \ar[dd] \ar[dl] \\
X \ar[dd]^f \ar[rrr] & & & \cU_h \ar[dd] \\
& \tilde{B} \ar[dl]^{\phi}_(0.5){\textrm{quasi-finite, surjective} \rightarrow} \ar[rrr]_{\eta}^(0.4){\textrm{surjective}} & & & D \ar[dl]_{\psi}^{\textrm{quasi-finite}}  \\
B \ar[rrr]_{\nu} & & & \fM_h
},
\end{equation*}
where all ``vertical'' squares are Cartesian, $\cU_h$ is the universal family over $\fM_h$ and $D$ is a smooth, proper scheme \cite[Theorem 9.25]{Viehweg_Quasi_projective_moduli}. Since all vertical maps are smooth, all relative canonical sheaves are compatible with pullbacks. By Corollary \ref{cor:rel_canonical_sheaf_big}, $\omega_{Y/D}$ is big. Hence
\begin{multline*}
 \kappa(\omega_{X/B}) 
= \underbrace{\kappa(\xi^* \omega_{X/B})}_{\parbox{85pt}{\tiny Kodaira dimension does not change by pulling back}}
= \underbrace{\kappa(\omega_{\tilde{X}/\tilde{B}})}_{\omega_{\tilde{X}/\tilde{B}} \cong \xi^* \omega_{X/B}} 
= \underbrace{\kappa( \zeta^* \omega_{Y/D})}_{\omega_{\tilde{X}/\tilde{B}} \cong \zeta^* \omega_{Y/D}} \\ = \kappa(\omega_{Y/D}) = \dim Y = \dim D +n = \var \tilde{f} + n = \var f +n 
\end{multline*}
\end{proof}

\section{Compact bases}
\label{sec:compact}

In this section the compact base case (i.e. Theorem \ref{thm:main_projective} and Corollary \ref{cor:main_projective}) is treated. Having a projective base allows us to use certain techniques not available in the general case. More precisely, the set of families of canonically polarized manifolds with fixed Hilbert polynomial form a nice moduli space if the base is projective. This is worded by the following lemma. However, first some preparation is necessary. 

Fix a projective manifold $B$, and a polynomial $h$. One can define a moduli functor $\fM_{B,h}$ of families of canonically polarized manifolds with Hilbert polynomial $h$ by
\begin{equation*}
\fM_{B,h}(T) : = \left\{ f: X \to B \times T \left| \parbox{180pt}{ $f$ is a smooth morphism, $\omega_f$ is $f$-ample, and $\chi(\omega_f^n|_{X_{(b,t)}})=h(n)$ for every $n \in \bZ$ and $(b,t) \in B \times T$  } \right. \right\}
\end{equation*}
One can also give a natural category fibered in groupoid structure to this functor which we also denote by $\fM_{B,h}$. Then the following lemma holds. 
 
\begin{lem}
\label{lem:hom_stack}
$\fM_{B,h} \cong \underline{\Hom}(B, \fM_h)$ as categories fibered in groupoids, where $\underline{\Hom}(B, \fM_h)$ is the $\Hom$-stack (\cite[Lines 1-4]{Olsson_Home_stacks_and_restriction_of_scalars}). In particular by \cite[Theorem 1.1]{Olsson_Home_stacks_and_restriction_of_scalars}, $\fM_{B,h}$ is a Deligne-Mumford stack, locally of finite type over $k$. 
\end{lem}

The next corollary is the reason why a locally of finite type DM stack structure on $\fM_{B,h}$ is useful.

\begin{cor}
If $f : X \to B$ is a family of canonically polarized manifolds with Hilbert polynomial $h$ over compact $B$, then $f$ is rigid (according to Definition \ref{def:rigid}) if  the infinitesimal deformation space $T^1(f,\fM_{B,h})$ of $f$ is zero.
\end{cor}

The expression for $T^1(f,\fM_{B,h})$ can be found for example in \cite[Theorem 1.1]{OMC_DT}. Then using that $\bL_{X/B} \cong \Omega_{X/B}$ in this case, one gets the following corollary.

\begin{cor}
\label{cor:inf_def_space_explicit}
A family $f : X \to B$ of canonically polarized manifolds over a compact base is rigid if $H^1(X, \sT_{X/B}) =0$.
\end{cor}

\begin{lem}
\label{lem:etale_trivialization}
If $f: Y \to S \times U$ is a family of canonically polarized manifolds with $U$ a manifold and $S$ a projective manifold, such that $W:=Y_u \to S \times \{u\} \cong S$ are isomorphic for all $u \in U$ as schemes over $S$, then there is a finite  \'etale cover $T \to U$ from a variety, such that the following isomorphisms holds
\begin{equation}
\label{eq:isomorphism_families}
\xymatrix{
Y \times_{S \times U} S \times T \cong Y \times_U T \ar[d] \ar@{<->}[r]^(0.7){\cong} & W \times T \ar[d] \\
S \times T \ar@{=}[r] & S \times T
}
\end{equation}
\end{lem}

\begin{proof}
Define the $S$-scheme $W:=Y_u \to S \times \{u\} \cong S$, for some choice of $u \in U$. Let $h$ be the Hilbert polynomial of $W \to S$. By Lemma \ref{lem:hom_stack} we know that $\fM_{S,h}$ is DM stack of finite type. In particular it has an \'etale (not necessarily finite) cover $\pi : V \to \fM_{S,h}$ by a scheme. Fix this cover $\pi$. Note that then $(\pi^{-1}([W \to S]))_{\red}$ is a zero dimensional reduced scheme of finite type, and then also proper over $\fM_{S,h}$. 

The family $Y \to S \times U$ defines a map $U \to \fM_{S,h}$ with zero dimensional image. Define $\tilde{U}:= U \times_{\fM_{S,h}} V$. Note that $\tilde{U}$ is a scheme, and $\tilde{U} \to U$ is \'etale. Hence, $\tilde{U}$ is smooth,  and then $\tilde{U} \to \fM_{S,h}$ factorizes through $(\pi^{-1}([W \to S]))_{\red}$. That is, in the definition of $\tilde{U}$, $V$ can be replaced by $(\pi^{-1}([W \to S]))_{\red}$. With other words, $\tilde{U} \cong U \times_{\fM_{h,S}} (\pi^{-1}([W \to S]))_{\red}$. In particular, $\tilde{U}$ is the disjoint union of  manifolds that are proper over $U$. Choose any of these, and define it to be $T$. Then $T$ factorizes through a point of $(\pi^{-1}([W \to S]))_{\red}$, which implies, that the associated family to $T \to \fM_{S,h}$ is the trivial family $W \times T \to S \times T$. However, $T \to \fM_{S,h}$ also factorizes through $U \to \fM_{S,h}$, which gives us the isomorphism \eqref{eq:isomorphism_families}.
\end{proof}

\begin{lem}
\label{lem:etale_trivialization2}
If $f: Y \to S \times U$ is a family of curves of genus at least two with $U$ a manifold and $S$ a projective manifold, such that for some $u \in U$, the restriction $W:=Y_u \to S \times \{u\}$ is non-isotrivial, then there is a finite  \'etale cover $T \to U$ from a variety, such that the following isomorphisms holds
\begin{equation*}
\xymatrix{
Y \times_{S \times U} S \times T \cong Y \times_U T \ar[d] \ar@{<->}[r]^(0.7){\cong} & W \times T \ar[d] \\
S \times T \ar@{=}[r] & S \times T
}
\end{equation*}
\end{lem}

\begin{proof}
%
First, notice that 
\begin{equation*}
H^1(W,\sT_{W/S} ) \cong H^1(W,\omega_{W/S}^{-1} ) =0,
\end{equation*}
by \cite[Corollary 5.12.c]{Esnault_Viehweg_Lectures_on_vanishing_theorems} and Corollary \ref{cor:rel_canonical_Kodaira_dim}. Hence, by Corollary \ref{cor:inf_def_space_explicit}, $W \to S$ is rigid. That is, for any $u \in U$, $Y|_{S \times \{u\}} \cong W$ as schemes over $S$. Applying Lemma \ref{lem:etale_trivialization} concludes the proof.
\end{proof}

\begin{lem}
\label{lem:stein_factorization}
If $f : X \to U$ is a smooth map onto a smooth curve, then for the  Stein-factorization
\begin{equation*}
\xymatrix{
 X \ar[r]^g \ar@/^2pc/[rr]^f & U' \ar[r]^h & U ,
}
\end{equation*}
the following holds: 
$U'$ is a smooth curve, $g$ is smooth and $h$ is \'etale.
\end{lem}

\begin{proof}
Since $X$ is normal, so is $U'$.
Then, one obtains smoothness of $U'$ by the equivalence of normality and smoothness in dimension one. For the rest of the statements, take any point $P \in X$. Then there is a diagram of tangent maps
\begin{equation*}
\xymatrix{
 T_{X,P} \ar[rr]^{T_{g,P}} \ar@/^2pc/[rrrr]^{T_{f,P}} & & T_{U',g(P)} \ar[rr]^{T_{h,g(P)}} & & T_{U,f(P)}
}.
\end{equation*}
Since the two tangent spaces on the right are one dimensional and the composition map is surjective, the only way to make the diagram commutative, if $T_{h,g(P)}$ is isomorphism  and $T_{g,P}$ is surjective. This proves everything stated in the proposition.
\end{proof}

\begin{proof}[Proof of Theorem \ref{thm:main_projective}]
First we prove, that if all $\var f_i \geq 1$, then $f$ is rigid. By Corollary \ref{cor:inf_def_space_explicit}, it is enough to show that $H^1(X,\sT_{X/U})=0$. Call $g_i$ the maps $f_{i+1} \circ \dots \circ f_n : X \to X_i$. Then $\sT_{X/U}$ has a filtration the quotient sheaves of which are $g_i^* \sT_{X_i/X_{i+1}}$. Hence it is enough to prove that $H^1(X,g_i^* \sT_{X_i/X_{i+1}})=0$ for all $i$. This follows from Proposition \ref{prop:rel_canonical_nef}, Corollary \ref{cor:rel_canonical_Kodaira_dim} and the vanishing theorem \cite[Corollary 5.12.c]{Esnault_Viehweg_Lectures_on_vanishing_theorems}.

We prove the other direction (or other statement) by induction on $n$. For $n=1$ it is true by Lemma \ref{lem:etale_trivialization}. So, assume that $n$ is arbitrary, and the statement is true for $n-1$. Then there are two possibilities. $f_n$ is either isotrivial or not. If it is isotrivial, then let $F$ be its fiber. By applying first Lemma \ref{lem:etale_trivialization} one obtains the upper Cartesian square of the following diagram, and then Lemma \ref{lem:stein_factorization} gives the lower factorization $X_{n-1}' \to U' \to U$.
\begin{equation*}
\xymatrix{
X \ar[d] & F \times X_{n-1}' \ar[d] \ar[l] \\
X_{n-1} \ar[d] & X_{n-1}' \ar[l]^{\textrm{\'etale}} \ar[d]^{\textrm{family of canonically polarized manifolds}} \\
U & \ar[l]^{\textrm{\'etale, finite}} U' 
}
\end{equation*}
As it is indicated on the diagram, $X_{n-1}' \to U'$ has canonically polarized fibers, since all its fibers are \'etale covers of fibers of $X_{n-1} \to U$. So, by setting $W:=X_{n-1}'$ and $Y:=F$, the inductional step is proven if $f_n$ is isotrivial. 

If $f_n$ is not isotrivial, then there must be some other $i$, for which $f_i$ is isotrivial. However, then using the inductional hypothesis, there is a diagram as follows.
\begin{equation*}
\xymatrix{
X_{n-1} \ar[d] & X_{n-1}':=W_{n-1} \times Y_{n-1} \ar[l]^(0.75){\textrm{\'etale}} \ar[d] \\
U & \ar[l]^{\textrm{\'etale, finite}} U''
}
\end{equation*}
Here $Y_{n-1}$ is a positive dimensional canonically polarized manifold, the map $W_{n-1} \to U''$ is a family of canonically polarized manifolds and $W_{n-1} \times Y_{n-1} \to U''$ is the composition of the first projection with the map $W_{n-1} \to U''$. Define $X':= X \times_{X_{n-1}} X_{n-1}'$. Since $X \to X_{n-1}$ is not isotrivial, same holds for $X' \to X_{n-1}'$. Then, there is either a $w \in W_{n-1}$ or a $y \in Y_{n-1}$, such that $X' \to X_{n-1}'$ is non-isotrivial over $\{w \} \times Y_{n-1}$ or $W_{n-1} \times \{y\}$. Assume first, that the first case is happening. Define then $Y:=X'|_{\{w \} \times Y_{n-1}}$. By Lemma \ref{lem:etale_trivialization2}, we obtain a Cartesian diagram as follows.
\begin{equation*}
\xymatrix{
X_{n-1} \ar[d] & &  \ar[ll]^{\textrm{\'etale, finite}} \ar[d] W \times Y  \\
X_{n-1}':=W_{n-1} \times Y_{n-1} & & \ar[ll]^(0.4){\textrm{\'etale, finite}} W \times Y_{n-1}
}
\end{equation*}
Then by taking Stein factorization of $W \to U''$  and using Lemma \ref{lem:stein_factorization}, we get the following diagram, where the preceding construction is also included and which proves the inductional step if $X' \to X_{n-1}'$ was non-isotrivial over $\{w \} \times Y_{n-1}$. 
\begin{equation*}
\xymatrix{
X \ar[d] & \ar[l]^{\textrm{\'etale, finite}} \ar[d] X' & W \times Y \ar[d] \ar[l]^{\textrm{\'etale, finite}} \\
X_{n-1} \ar[d] & \ar[l] X_{n-1}'=W_{n-1} \times Y_{n-1} \ar[d] & \ar[d] W \times Y_{n-1} \ar[l]^(0.3){\parbox{25pt}{\tiny \'etale, finite}} \ar[d] \\
U & \ar[l]^{\textrm{\'etale, finite}} U'' & \ar[l]^{\textrm{\'etale, finite}} U'
}
\end{equation*}
We conclude the proof with the case when $X' \to X_{n-1}'$ is non-isotrivial over $W_{n-1} \times \{y\}$. Then, define $W:= X'|_{W_{n-1} \times \{y\}}$. Using Lemma \ref{lem:etale_trivialization2} again yields the following diagram, where the left top square is Cartesian.
\begin{equation*}
\xymatrix{
X \ar[d] & \ar[l] \ar[d] X' & W \times Y_{n-1}' \ar[d] \ar[l] \\
X_{n-1} \ar[d] & \ar[l] X_{n-1}'=W_{n-1} \times Y_{n-1} \ar[d] &  W_{n-1} \times Y_{n-1}' \ar[l]  \\
U & \ar[l] U''
} 
\end{equation*}
Setting $Y:=Y_{n-1}'$ and $U':=U''$ yields the result in this case too.
\end{proof}

\begin{proof}[Proof of Corollary \ref{cor:main_projective}]
Corollary \ref{cor:main_projective} is really a corollary of the Proof of Theorem \ref{thm:main_projective} instead of the statement itself. Recall that we really proved not only the rigidity of $U \to \fM_h$ in the situation of Notation \ref{notation:tower} when $\var f_i \geq 1$, but also its infinitesimal rigidity. 

Let us get back now to the statement of Corollary \ref{cor:main_projective}. 
First we show that  there is a  lifting as follows, where $U'$ is another smooth curve mapping finitely to $U$. 
\begin{equation}
\label{eq:main_projective:lifting}
\xymatrix{
U' \ar[r]^{\nu'} \ar[d]^{\tau} & \widetilde{\fK \fF}_h \ar[d] \\
U \ar[r]^{\nu} & \fM_h ,
}
\end{equation}
Indeed, by Theorem \ref{thm:proper}, $\widetilde{\fK \fF}_h \to \fM_h$ is proper and by Lemma \ref{lem:representable}, it is also representable by algebraic spaces. Hence if we set $U':= U \times_{\fM_h} \widetilde{\fK \fF}_h$, then $U'$ is an algebraic space, finite and \'etale over $U$. In particular, it is smooth of dimension one, and hence it is a scheme, more precisely a projective curve.

Given a lifting as in \eqref{eq:main_projective:lifting}, $\nu \circ \tau$ is infinitesimally rigid according to the proof of Theorem \ref{thm:main_projective}. That is, for the induced family $f': X' \to U'$, $H^1(X',\sT_{X'/U'})=0$. Let $f : X \to U$ be the family induced by $\nu$. Then $X' = X \times_U U'$, so in particular $\sT_{X'/U'} \cong (\tau')^* \sT_{X/U}$, where $\tau'$ is the natural projection map $X' \to X$. So, the following computation implies that $\nu$ is infinitesimally rigid and hence rigid as well.
\begin{multline*}
0 = H^1( X', \sT_{X'/U'}) 
=\underbrace{ H^1( X, \tau'_* \sT_{X'/U'} )}_{\textrm{  $\tau$ and hence $\tau'$ is finite}}
=\underbrace{ H^1( X, \tau'_* (\tau')^* \sT_{X/U} )  }_{ \sT_{X'/U'} \cong (\tau')^* \sT_{X/U}}
\\ =\underbrace{H^1( X, \tau'_* \sO_{X'} \otimes \sT_{X/U} )}_{\textrm{  projection formula}}
 \hookleftarrow \underbrace{H^1( X, \sO_X \otimes \sT_{X/U} )}_{\parbox{100pt}{ \tiny splitting $\sO_X \to \tau'_* \sO_{X'}$ given by the trace map}}
 = H^1 (X, \sT_{X/U})
\end{multline*}
\end{proof}

\section{Relative properness of the moduli space of iterated Kodaira fibrations}
\label{sec:relative_properness}

In this section we show that the natural forgetful map $\pi : \widetilde{\fK\fF_h} \to \fM_h$, defined in Definition \ref{def:iterated_Kodaira_Spencer}, is proper. In particular, $\fK\fF_h$ is closed in $\fM_h$, and then consequently it is a connected component of $\fM_h$.  This is used for example to conclude Corollary \ref{cor:main_projective} from Theorem \ref{thm:main_projective}.

\begin{notation}
\label{notation:iteated_Kodaira_Fibration_compactification}
Choose a connected component $\fY$ of $\widetilde{\fK\fF_h}$. Consider a tower
\begin{equation*}
\xymatrix{
X= X_n \ar[r]^{f_n} & X_{n-1} \ar[r]^{f_{n-1}} & \dots \ar[r]^{f_2} & X_1 \ar[r]^-{f_1} & X_0= \Spec k \in \fY(\Spec k),
}
\end{equation*}
and the following two groups of numerical invariants: the genus $g_i$ of the fibers of $f_i$ and the degree $d_i:=\deg \left( \nu_i|_{(X_{i-1})_y} \right)^* \lambda_{g_i}$, where $y \in X_{i-2}$ is arbitrary, $\nu_i : X_{i-1} \to \ofM_{g_i}$ is the moduli map of $f_i$, and $\lambda_g$ is an ample line bundle on $\ofM_g$. By flatness these invariants are independent of the chosen tower. Hence, there is a proper DM-stack 
\begin{equation}
\label{eq:compactified_moduli}
\fM:=\fM_{g_1, 0}^{\textrm{balanced}}( \fM_{g_2,0}^{\textrm{balanced}}( \dots ( \ofM_{g_n,0}, d_n) \dots ), d_2), 
\end{equation}
compactifying $\fY$, given by the iterated use of the Abramovich-Vistoli construction of stable maps 
\cite{Abramovich_Vistoli_Complete_moduli_for_fibered_surfaces,Abramovich_Vistoli_Compactifying_the_space_of_stable_maps}. Let $\ofY$ be the closure of $\fY$ in $\fM$. Note that we used the characteristic zero assumption for $\fM$ to be a DM-stack.
\end{notation}

The objects of $\fM$ defined in Notation \ref{notation:iteated_Kodaira_Fibration_compactification} over a test scheme $U$ are towers
\begin{equation}
\label{eq:stacky_tower}
\xymatrix{
\sX= \sX_n \ar[r]^{\tilde{f}_n} & \sX_{n-1} \ar[r]^{\tilde{f}_{n-1}} & \dots \ar[r]^{\tilde{f}_2} & \sX_1 \ar[r]^-{\tilde{f}_1} & \sX_0= U
},
\end{equation}
where 
\begin{enumerate}
\item all $\sX_i$ are  DM-stacks,
\item $\tilde{f}_n$ is a family of stable curves
\item for $1 \leq i \leq n-1$, $\tilde{f}_i$ are  flat families of DM-stacks with one dimensional, geometrically connected fibers with only smooth and nodal singularities,
\item for $1 \leq i \leq n-1$, the coarse moduli spaces of the fibers of $\tilde{f}_i$ are stable curves and the only possible stack structure of these fibers  appears over the nodes, where it is \'etale locally isomorphic to a (stack) quotient of $  \Spec \left( \factor{ k[x,y]}{(xy)} \right)  $, by the action of $\mu_r$  via $(\xi x, \xi^{-1} y)$, and finally
\item  the induced $\mu_r$ action on the fibers of
\begin{equation*}
\xymatrix{
\sX= \sX_n \ar[r]^{\tilde{f}_n} & \sX_{n-1} \ar[r]^{\tilde{f}_{n-1}} & \dots \ar[r]^{\tilde{f}_{i+1}} & \sX_i 
},
\end{equation*}
as points of 
\begin{equation*}
\fM_{g_{i+1}, 0}^{\textrm{balanced}}( \fM_{g_{i+2},0}^{\textrm{balanced}}( \dots ( \ofM_{g_n,0}, d_n) \dots ), d_{i+2}), 
\end{equation*}
over the nodes of the fibers of $f_i$ has to be essential, that is, faithful. 
\end{enumerate}
It should be noted that the last condition is included only for the sake of completeness and will not be used in our proofs.

One may also consider the coarse moduli tower of a tower as in \eqref{eq:stacky_tower}, obtaining a commutative diagram as follows.
\begin{equation}
\label{eq:schemey_tower}
\xymatrix{
\sX= \sX_n \ar[r]^{\tilde{f}_n} \ar[d]^{ \pi_n} & \sX_{n-1} \ar[r]^{\tilde{f}_{n-1}} \ar[d]^{\pi_{n-1}} & \dots \ar[r]^{\tilde{f}_2} & \sX_1  \ar[r]^-{\tilde{f}_1} \ar[d]^{\pi_{1}} & \sX_0= U \ar@{=}[d] \\
X= X_n \ar[r]^{f_n} & X_{n-1} \ar[r]^{f_{n-1}} & \dots \ar[r]^{f_2} & X_1 \ar[r]^-{f_1} & X_0= U
},
\end{equation}
 This process  would yield a morphism $\fM \to \ofM_h$ to the moduli space of stable varieties with Hilbert polynomial $h$, as soon as we proved that $X \to U$ of \eqref{eq:schemey_tower} is a stable family.   This is shown in Theorem \ref{thm:stable}. For a quick introduction to stable varieties we refer to \cite{Kollar_Moduli_of_varieties_of_general_type} and \cite{Kollar_Singularities_of_the_minimal_model_program}. 

\begin{notation}
\label{notation:iteated_Kodaira_Fibration_compactified_map}
In the situation of Notation \ref{notation:iteated_Kodaira_Fibration_compactification}, Let $\overline{\rho} : \fM \to \ofM_h$ be the morphism given by the above considerations and $\rho:= \overline{\rho}|_{\fY}$. By abuse of notation we also denote by $\overline{\rho}$ the restriction $\overline{\rho}|_{\overline{\fY}}$.
\end{notation}

Before proving Theorem \ref{thm:stable}, we collect all the  properties of the towers \eqref{eq:stacky_tower} that are used in the following proofs.

\begin{fact}
\label{fact:towers}
For any scheme $U$ over $k$ and any element of $\fM(U)$, using the notations of \eqref{eq:schemey_tower}, the following holds.
\begin{enumerate}
\item \label{itm:towers:smooth_towers} For the \emph{smooth towers}, i.e., for which all  $\tf_i$ are smooth, $\sX_i \to X_i$ is an isomorphism for all $i$. 
\item \label{itm:towers:isomorphism_on_the_first_level} If $X_1$ is smooth over $U$ at a point $P$, then $\sX_1 \to X_1$ is isomorphism over $P$.
\item \label{itm:towers:general_points} For each $u \in U$, the coarse moduli map $\sX_u \to X_u$ is an isomorphism at the general points of $\sX_u$ and furthermore, these points are smooth. 
\item \label{itm:towers:codimension_one_points_singular} For each $u \in U$, the singular (i.e., not regular) codimension one points of $\sX_u$ map onto singular codimension one points of $X_u$.
\item \label{itm:towers:codimension_one_points_regular} For each $u \in U$, the map $\sX_u \to X_u$ is an isomorphism over the regular codimension 1 points  of $X_u$. 
\end{enumerate}

\end{fact}

\begin{proof}
\begin{enumerate}
\item Follows from the definition.
\item By definition, the points of $\sX_1$  where $\sX_1 \to X_1$ is not an isomorphism as well as the images of these points in $X_1$ are nodal over $U$.
\item Immediate from the fact that $\tf_i$ are smooth families of schemes in relative codimension zero over $\sX_{i-1}$. Note also that taking coarse moduli space is compatible by base-change according to \cite[Corollary 3.3]{Abramovich_Olsson_Vistoli_Tame_stacks_in_positive-characteristic}.
\item  Instead of writing $u$ in subindex everywhere in the proof, let us assume that $U$ is the spectrum of a field (of characteristic zero). We show the statement by induction on $n$. If $n=1$, the statement is automatic, since $\sX = X$. So, assume $n>1$. Let $P$ be a codimension one point of $\sX$, and $Q$ its image in $\sX_1$. There are two cases. 
\begin{itemize}
\item If $Q$ is a smooth point, then $\sX_1 \to X_1$ is isomorphism at $Q$ (by point \eqref{itm:towers:isomorphism_on_the_first_level})  and $P$ is either a general or a codimension 1 point in $\sX_Q$. In the former case $P$ is smooth in $\sX_Q$ by point \eqref{itm:towers:general_points}, and then it is also smooth in $\sX$, which cannot happen. On the other hand, in the latter case  $P$ has to be not only singular in $X$ but also in $X_Q$.  Hence, it maps to a singular point of $X_Q$ by the inductional assumption. However, this means that $P$ maps to a singular point of $X$ as well, because $Q$ has to be a general point of $X_1$ and hence $\sO_{X_Q,P} \cong \sO_{X,P}$.
\item If $Q$ is a node. Then $P$ is a general point of $\sX_Q$. However, then the image $Q'$ of $Q$ in $X_1$ is also a node, and $P$ maps to a general point $P'$ of $X_{Q'}$. So, $X$ is flat over a node at $P'$, which implies that it has to be singular at $P'$. 
\end{itemize}
\item Again, instead of writing $u$ in subindex everywhere in the proof, let us assume that $U$ is the spectrum of a field. We also prove by induction on $n$. From the analysis of the previous point we see that if $P' \in X$ is a regular codimension 1 point, then its image $Q$ in $X_1$ has to be regular as well, and furthermore, $\sX_1 \to X_1$ is isomorphism at $Q$. There are two cases depending on whether $Q$ has codimension $1$ or $0$. In the first case $P'$ is a general point in $X_{Q}$, and hence by point \eqref{itm:towers:general_points} $\sX_{Q} \to X_{Q}$ is isomorphism over $P'$.  In the second case $P'$ is a codimension one point of $X_{Q}$. Since $X_{Q}$ is a general fiber and $P'$ is regular in $X$, it also has to be regular in $X_{Q}$. However, then by the inductional hypothesis, we obtain the same conclusion: $\sX_{Q} \to X_{Q}$ is isomorphism over $P'$. So, in both cases we have arrived to this same conclusion: $\sX_{Q} \to X_{Q}$ is isomorphism over $P'$. We claim that using that $P'$ is a 
regular point,  it follows that also $\sX \to X$ is 
an isomorphism over $P'$. Using \cite[Lemma 2.2.3]{Abramovich_Vistoli_Compactifying_the_space_of_stable_maps}, this is equivalent to showing that if $Y \to Z$ is a morphism of schemes over the scheme $B$, $P \in Z$ is a regular point mapping to $Q$ in $B$, $Y$ is Cohen-Macaulay and $Y_Q \to Z_Q$ is \'etale over $P$, then so is $Y \to Z$. However, the above translation is true, since by the regularity and Cohen-Macaualayness assumption $Y \to Z$ and $Y_Q \to Z_Q$ are flat \cite[Exercise III.10.9]{Hartshorne_Algebraic_geometry}, and hence \'etaleness at $P$ is equivalent to unramifiedness \cite[III.10.3]{Hartshorne_Algebraic_geometry}. However, unramifiedness over  $P$ of $Y_Q \to Z_Q$ and $Y \to Z$ is equivalent. 
\end{enumerate}

\end{proof}

\begin{thm}
\label{thm:stable}
For any element of $\fM(U)$ as in \eqref{eq:stacky_tower}, the family $X \to U$ obtained by taking coarse moduli spaces as in \eqref{eq:schemey_tower} is a stable family.
\end{thm}

\begin{proof}
First, we review what a stable family is. A scheme $W$ is a stable scheme if it has slc  singularities \cite[Definition-Lemma 5.10]{Kollar_Singularities_of_the_minimal_model_program}  and $\omega_W$ is ample. Note that $\omega_W$ in this case is not a line bundle, only a $\bQ$-line bundle, that is $\omega_{W}^{[m]}:= (\omega_{W}^{m})^{**}$ is a line bundle for some $m \in \bZ_{>0}$. A flat family $V \to Y$ is stable, if all its fibers are stable schemes and $\omega_{V/Y}^{[m]}$ is flat and compatible with arbitrary base change for every $m \in \bZ$. 

After the above preliminary review, we may start the proof. Fix an element of $\fM(U)$ as in \eqref{eq:schemey_tower}. First, \emph{we claim that  $\omega_{X/U}^{[m]}$ is a flat and relatively Cohen-Macaulay $\bQ$-line-bundle for every $m \in \bZ$}. Note at this point that since $\sX \to U$ is a tower of relatively Gorenstein families, $\sX$ is relatively Gorenstein over $U$ as well. That is, it is relatively Cohen-Macaulay over $U$ and $\omega_{\sX/U}$ is a line bundle. Furthermore, by Corollary \ref{cor:pullback_K_X}, $\omega_{\sX/U}^{m} \cong \pi^{[*]} \omega_{X/U}^{[m]}$. In particular, then by Lemma \ref{lem:stack_reflexive_S_r}, using that the fibers of $X \to U$ are $G_1$ and $S_2$ because of Lemma \ref{lem:demi_normal}, $\omega_{X/U}^{[m]}$ is relatively Cohen-Macaulay and flat over $U$ for every $ m \in \bZ$. 
Note now that since $\sX$ is a DM-stack, there is an $N$, such that $\omega_{\sX/U}^N$ descends to a line bundle $\sL$ on $X$. That is, for this $N$,  $\pi^* \sL \cong \omega_{\sX/U}^N$. However then, since $\omega_{X/U}$ is a line bundle in relative codimension one and also in relative codimension one the isomorphism $\pi^* \omega_{X/U} \cong \omega_{\sX/U}$ holds, in codimension one we have  $\omega_{X/U}^{[N]} \cong \pi_* \omega_{\sX/U}^N \cong \sL$. Thus, since both $\omega_{X/U}^{[N]}$ and $\sL$ are flat and relatively $S_2$, they are isomorphic globally by \cite[Proposition 3.5]{Hassett_Kovacs_Reflexive_pull_backs}. That is, $\omega_{X/U}^{[N]}$ is a $\bQ$-line bundle indeed, which finishes the proof of the claim.  

Now, we use the claim to show that $X \to U$ is stable. First, to show that for every $u \in U$,  $X_u$ is slc, by Lemma \ref{lem:stack_demi_normal}, Lemma \ref{lem:demi_normal} and Lemma \ref{lem:finite_Galois_slc}.\ref{itm:finite_Galois_slc:slc}, we are supposed to show that $K_{X_u}$ is $\bQ$-Cartier.  However, this is the special case of the above claim when the base is a point.

Next step is to show that $\omega_{X_u}$ is ample for every $u \in U$. Since, by the proof of the claim, $\pi^* \omega_{X_u}^{[N]} \cong \omega_{\sX_u}^{N}$, it is enough to show that for any tower as in \eqref{eq:stacky_tower}, $\omega_{\sX_i/\sX_{i-1}}$ is nef  and $\tilde{f}_i$-ample for all $2 \leq i \leq n$ (recall that a line bundle is ample on a DM-stack if some of its power descends to an ample line bundle on the coarse moduli space). By \cite[Proposition 2.6]{VA_ITO}, there is a finite cover $Y \to \sX_{i-1}$ by a scheme. Then, $\omega_{(\sX_i)_Y/Y}$ is the pullback of $\omega_{\sX_i/\sX_{i-1}}$, hence it is enough to prove the nefness and relative ampleness of the former. However, that is obvious, since $(\sX_i)_Y \to Y$ is a family of stable curves over a scheme. This concludes the proof of the stability of $X_u$. 

We are left to show that $\omega_{X/U}^{[m]}$ is flat and compatible with base change. We have already proved that it is relatively Cohen-Macaualay and flat. Then, \cite[Corollary 3.8]{Hassett_Kovacs_Reflexive_pull_backs}  concludes our proof.
\end{proof}

The following is the main theorem of the section. See Definition \ref{defn:iterated_Kodaira_fibration} for the definition of the space of iterated Kodaira fibrations $\widetilde{\fK\fF}_h$ and of the forgetful map $\pi$.

\begin{thm}
\label{thm:proper}
The natural forgetful map $\pi : \widetilde{\fK\fF}_h \to \fM_h$ is proper. 
\end{thm}

\begin{proof}
It is enough to show that $\pi$ is proper on every component of $\widetilde{\fK\fF}_h$. 
So, choose a connected component $\fY$ of $\widetilde{\fK\fF}_h$ as in Notation \ref{notation:iteated_Kodaira_Fibration_compactification}. Using Notations \ref{notation:iteated_Kodaira_Fibration_compactification} and \ref{notation:iteated_Kodaira_Fibration_compactified_map}, $\pi|_{\overline{\fY}}$ agrees with the map $\overline{\rho}$, which according to  Theorem \ref{thm:stable}, can be put into  a commutative diagram of Deligne-Mumford stacks as follows.
\begin{equation*}
\xymatrix{
 \fY \ar[r]  \ar[d]^{\rho} & \ofY \ar[d]^{\overline{\rho}} \\
\fM_h \ar[r] & \ofM_h
} 
\end{equation*}
Since both $\ofY$ and $\ofM_h$ are proper, so is $\overline{\rho}$. In particular, then $\overline{\rho}^{-1}(\fM_h)$ is proper over $\fM_h$. Therefore, it is enough to show that the two open substacks $\fY$ and $\overline{\rho}^{-1}(\fM_h)$ of $\ofY$ coincide, or equivalently that the natural map between the isomorphism classes of their $k$-points is a bijection. 

Fix a tower as in \eqref{eq:stacky_tower} over $U=\Spec k$, such that
\begin{enumerate}
\item \label{itm:proper:deformation}  there is a one-parameter family 
\begin{equation*}
\xymatrix{
X' = X_n' \ar[r]^{f_n'}  & X_{n-1}' \ar[r]^{f_{n-1}'}  & \dots \ar[r] & X_1' \ar[r]^{f_1'}  & X_0' = V \ar@{=}[d] \\
\sX'= \sX_n' \ar[r]^{\tilde{f}_n'} \ar[u]^{\pi_n'} & \sX_{n-1}' \ar[r]^{\tilde{f}_{n-1}'} \ar[u]^{\pi_{n-1}'} & \dots \ar[r]  & \sX_1' \ar[r]^{\tilde{f}_1'} \ar[u]^{\pi_1'} & \sX_0'= V ,
}
\end{equation*}
the fiber of which over $0 \in V$ is our fixed tower and its coarse moduli space tower (as in \eqref{eq:schemey_tower}) and over any other $v \in V$ the fibers are smooth towers (as defined in Fact \ref{fact:towers}.\ref{itm:towers:smooth_towers}, and 
\item furthermore the total space $X=X'_0$ of  is smooth.
\end{enumerate}
\emph{We are supposed to prove that the central fiber is a smooth tower as well, that is, all the morphisms $\tilde{f}_i$ of the central stacky tower are smooth.} 

We prove this statement by induction. For $n=1$ the above statement is obvious, since then smoothness of $X$ implies that $\sX = X$. For $n >1$, there is a smooth dense open subset $W \subseteq X_1$, such that $X \to X_1$ is smooth over $W$ (here we use the notations of \eqref{eq:schemey_tower}). Furthermore by Fact \ref{fact:towers}.\ref{itm:towers:isomorphism_on_the_first_level},    $\sX_1 \to X_1$ is an isomorphism  over $W$. Fix then a closed point $w \in W$ and choose a smooth curve $ C \to X_1'$ through $w$ meeting $X_1$ transversally. Pulling back everything over $C$, we put ourselves in the same situation but for $n$ replaced by $n-1$. Hence, by the inductional hypothesis all $(\tf_i)_w$ are smooth. Hence by the openness of smoothness, $\tf_i$ are smooth over $W$. In particular for $i \geq 2$, $\tf_i : \sX_i \to \sX_{i-1}$ are smooth over a non-empty open set of $\sX_{i-1}$. Furthermore, $X_i$ are irreducible, since they are 
the surjective image of the irreducible $X$, and then $\sX_i$ are irreducible as well. Therefore, Lemma \ref{lem:degeneration_Kodaira_fibration} applies, whence $\tf_i$ are smooth for $i \geq 2$. But then if there was a nodal point of $\sX_1$ (or equivalently of $X_1$), $X$ could not be smooth, which finishes our proof.
\end{proof}

\begin{rem}
In Definition \ref{defn:iterated_Kodaira_fibration}, in the definition of iterated Kodaira fibrations $\Var \left( f_i|_{X_u} \right)$ is assumed. This assumption is not used anywhere in the current section and also not in \cite[Theorem 1.2]{Patakfalvi_Fibered_stable_varieties}, where it is shown that $\widetilde{\fK\fF}_h \to \fK\fF_h$ is \'etale. Hence the results of the current section and the consequence that $\fK\fF_h$ consists of connected components of $\fM_h$ are valid without assuming the above variational assumption.
\end{rem}

\subsection{Lemmas}

Here we list a few lemmas used in the main statements of this section. The first one is specifically designed to prove that $\omega_{X/U}^{[m]}$ is relatively Cohen-Macaulay for coarse moduli towers as in \eqref{eq:schemey_tower}.

\begin{lem}
\label{lem:stack_reflexive_S_r}
Consider the following setup.
\begin{itemize}
\item $\sX$ is a (tame) Deligne-Mumford stack, flat over a scheme $U$, with coarse moduli space $\pi : \sX \to X$,
\item $X$ is relatively $S_2$ and $G_1$,
\item $\sF$ is a reflexive coherent sheaf on $X$, locally free in relative codimension one,
\item the reflexive pullback $\pi^{[*]} \sF:=(\pi^* \sF)^{**}$ is flat and relatively $S_r$  for some $r \geq 2$.
\end{itemize}
In the above situation $\sF$ is a flat and relative $S_r$ sheaf over $U$.
\end{lem}

\begin{rem}
Note that a coherent sheaf by definition is $S_r$ on a Deligne-Mumford stack if it is $S_r$ on some (or equivalently every) \'etale atlas. The property $G_1$ (Gorenstein in codimension one) is defined the same way. The word relative means that instead of requiring $\sF$ be $S_r$ or $G_1$ on the total space, we require it to be $S_r$ or $G_1$ on every fiber. 
\end{rem}

\begin{proof}[Proof of Lemma \ref{lem:stack_reflexive_S_r}]
There is a natural map $\alpha : \pi^* \sF \to (\pi^* \sF)^{**}$. The composition of the pushforward of this map with the other natural map $ \sF \to \pi_* \pi^* \sF$ gives a homomorphism $ \beta : \sF \to \pi_* (\pi^{[*]} \sF) $. Note that since $\sF$ is locally free in relative codimension one, the natural homomorphism $\alpha$ is isomorphism in relative codimension one. By the projection formula, using that $\pi_* \sO_{\sX} \cong \sO_X$, the same holds for $\sF \to \pi_* \pi^* \sF$, and hence for $\beta$. Furthermore, by \cite[Lemma 3.7]{Patakfalvi_Fibered_stable_varieties} and our assumption, $\pi_* (\pi^{[*]} \sF ) $ is flat and relatively $S_r$. That is, $\beta$ is a morphism between a reflexive and a flat, relative $S_2$ sheaf which is isomorphism in relative codimension one. In particular, then it is an isomorphism by \cite[Proposition 3.6.2]{Hassett_Kovacs_Reflexive_pull_backs}, which shows that $\sF$ is flat and relatively $S_r$  indeed. 
\end{proof}

The following lemmas prove demi-normality of certain schemes and stacks. The role of demi-normal for slc singularities is the same as of normal for log canonical singularities. The difference between the two is that demi-normal allows slightly worse singularities in codimension one, and consequently for example multiple irreducible components. A scheme $X$ is \emph{demi-normal} if it is $S_2$ and all its points of codimension at most one are regular or nodal \cite[Definition 5.1]{Kollar_Singularities_of_the_minimal_model_program}. Note also that if $U \to V$ is a surjective \'etale map, then $U$ is demi-normal if and only if so is $V$. Hence we may define a stack to be demi-normal if  it has a  demi-normal \'etale atlas. In this case all \'etale neighborhoods of the stack will be demi-normal. 

\begin{lem}
\label{lem:stack_demi_normal}
The towers as in \eqref{eq:stacky_tower} are demi-normal.
\end{lem}

\begin{proof}
It follows by induction on $n$ using that  $\tf_i$ are families of nodal curves.
\end{proof}

The corresponding statement for the coarse towers as in \eqref{eq:schemey_tower} is slightly more involved and it also uses the following two lemmas.

\begin{lem} \cite[Claim 1.41.1]{Kollar_Singularities_of_the_minimal_model_program}
\label{lem:regular_or_nodal}
Let $(A, m)$ be a local ring of dimension 1 with residue field $K$ and normalization $\oA$. Further, assume that $A$ is the quotient of a regular local ring. Then $A$ is nodal or regular if and only if $m$ is the intersection of all maximal ideal in $\oA$ over $m$, and $\dim_K \left( \factor{\oA}{m} \right) \leq 2$. 
\end{lem}

\begin{lem}
\label{lem:quotient_demi_normal}
If $X$ is demi-normal scheme of finite type over $k$ (note $k$ is assumed to be of characteristic 0) with an action of a finite group $G$, then $X/G$ is demi-normal as well.
\end{lem}

\begin{proof}
Set $Y:= X/G$. Then $Y$ is $S_2$ by \cite[Proposition  5.4]{Kollar_Mori_Birational_geometry_of_algebraic_varieties}.  We are supposed to show that for every codimension one point $Q \in Y$, $Y$ is either regular or nodal at $Q$ (nodal here is meant in the sense of \cite[1.41]{Kollar_Singularities_of_the_minimal_model_program}). So, fix such a point $Q$. Let $P$ be any of its preimages, and $H$ the stabilizer of $P$ in $G$. Then $X/H \to Y$ is \'etale over $Q$. Hence, we may replace $G$ by $H$, $Y$ by $X/H$ and $Q$ by any of the preimages of $Q$ in $X/H$. Equivalently, we may assume that $P$ is a fixed point of $G$. By dropping the finite typeness assumption, we may further assume that $X = \Spec A$ for a local ring $(A,m)$ of dimension 1, and that $P$ is the unique closed point of $X$. Let $\oA$ be the normalization of $A$, and $K$ the residue field of $A$. Note, that  $Y$ is the spectrum of the local ring $(A^G, m^G)$, the normalization of $Y$ is $\Spec \oA^G$ and $K^G$ is the residue field of $A^G$.  
Consider the following isomorphisms 
over $K^G$.
\begin{equation*}
 \factor{A^G}{m^G} \cong \factor{A^G}{m \cap A^G} 
\underbrace{\cong \factor{A^G + m }{m}}_{\parbox{50pt}{\tiny second isomorphism theorem}} 
\underbrace{\cong \left( \factor{A}{m} \right)^G}_{\parbox{70pt}{\tiny the action of $G$ is completely reducible}}
\end{equation*}
This implies the inequality
\begin{equation*}
\dim_{K^G} \left( \factor{A^G}{m^G} \right) = \dim_{K^G} \left( \factor{A}{m} \right)^G \leq \dim_K \left( \factor{A}{m} \right) .
\end{equation*}
Hence by Lemma \ref{lem:regular_or_nodal}, we are left to show that $m^G$ is the intersection of all maximal ideals lying above $m^G$ in $\oA^G$. Indeed, let $\{m_i\}$ be the set of maximal ideals over $m$ in $\oA$. Then $\left\{m_i \cap \oA^G \right\}$ is the set of maximal ideals lying over $m^G$ in $\oA^G$. Therefore, the intersection of the maximal ideals of $\oA^G$ lying over $m^G$ is
\begin{equation*}
\bigcap \left( m_i \cap \oA^G \right) =  \left( \bigcap m_i  \right) \cap \oA^G = m  \cap \oA^G = m \cap A^G = m^G.
\end{equation*}
Here, the third inequality follows from $m \subseteq A$ and its consequence $m \cap \oA^G \subseteq A \cap \oA^G= A^G$.

\end{proof}

\begin{lem}
\label{lem:demi_normal}
The formulation of the coarse moduli spaces as in \eqref{eq:stacky_tower} and \eqref{eq:schemey_tower} are compatible with base change (of $U$), conserves flatness (over $U$) and the coarse fibers (over $U$) are demi-normal. 
\end{lem}

\begin{proof}[Proof of Lemma \ref{lem:demi_normal}]
The compatibility with base change and conservation of flatness is shown in \cite[Corollary 3.3]{Abramovich_Olsson_Vistoli_Tame_stacks_in_positive-characteristic}. Hence, in particular, the fibers of the coarse families are the coarse moduli spaces of the fibers. So, we may assume that $U$ is the spectrum of a field. 
%

Next, note that since \'etale locally every fibers of $X$ is a quotient of an \'etale chart of $\sX$ by a finite group \cite[Lemma 2.2.3]{Abramovich_Vistoli_Compactifying_the_space_of_stable_maps}, \'etale locally every fiber is a quotient of a demi-normal scheme of finite type over $k$ by a finite group. In particular, it is demi-normal by Lemma \ref{lem:demi_normal}. 

\end{proof}

The next lemma is used to show that from $\sX$ of \eqref{eq:schemey_tower} being slc follows that $X$ is slc as well. Slc here means semi-log canonical \cite[Definition-Lemma 5.10]{Kollar_Singularities_of_the_minimal_model_program}.

\begin{lem}
\label{lem:finite_Galois_slc}
If $f: X \to Y$ is a finite, Galois cover of equidimensional, demi-normal schemes, then 
\begin{enumerate}
\item \label{itm:finite_Galois_slc:slc} if $K_X$ and $K_Y$ are $\bQ$-Cartier and $X$ is slc, then so is $Y$, and
\item \label{itm:finite_Galois_slc:unramified} if $f$ is unramified at the smooth codimension one points and there is no nodal codimension one point of $X$ mapping onto a regular codimension one point of $Y$, then for any compatible choice of canonical divisors, $f^* K_Y = K_X$ (pulling back Weil-divisors via finite maps make sense 
\cite[proof of Prop. 5.20]{Kollar_Mori_Birational_geometry_of_algebraic_varieties} so we do not need any $\bQ$-Cartier assumption here).
\end{enumerate}

\end{lem}

\begin{proof}
First, we set up the problem, and relate $f^* K_Y$ to $K_X$, a computation that will be used in the proof of both points. 
Let $\pi : (\oX,D) \to X$ and $\rho : (\oY,E) \to Y$ be the normalizations. Hence, there is a commutative diagram as follows
\begin{equation*}
\xymatrix{
X \ar[d]^f & \ar[l]^{\pi} (\oX,D) \ar[d]^{\of} \\
Y  & \ar[l]^{\sigma} (\oY,E) \\
},
\end{equation*}
where $\of$ is also a Galois cover.  Let $\Delta_i$ be the prime divisors over which $\of$ is ramified and let $r_i$ be the ramification index of $\Delta_i$. Then 
\begin{equation}
\label{eq:finite_Galois_slc:canonicals}
K_{\oX}  = \of^* ( K_{\oY} + \sum_{i} \frac{r_i -1 }{r_i} \Delta_i). 
\end{equation}
Consider now two codimension one points $x \in X$ and $y \in Y$ such that $f(x) = y$. Note that if $x$ is regular point, then there is a normal neighborhood of $f^{-1}(y)$. Hence the quotient $Y$ is also normal in a neighborhood of $y$, so $y$ is also a regular point. On the other hand if $x$ is nodal, then $y$ can be either regular or nodal. In particular, $D \geq D'$ for $D':= (\of^*E)_{red}$. Define then $E'$ to be the sum of the prime divisors appearing in $E$ over which $\of$ is not ramified. Then,  
\begin{equation}
\label{eq:finite_Galois_slc:boundaries}
D' = \of^* \left(E' + \sum_{\Delta_i \leq E - E'} \frac{1}{r_i} \Delta_i \right) .
\end{equation} 
Hence, 
\begin{equation}
\label{eq:finite_Galois_slc:computation}
\begin{split}
K_{\oX} + D' &  
=\underbrace{  \of^* \left( K_{\oY} + E' + \sum_{\Delta_i \leq E - E'} \frac{1}{r_i} \Delta_i + \sum_{i} \frac{r_i -1 }{r_i} \Delta_i \right)}_{
\textrm{ by \eqref{eq:finite_Galois_slc:canonicals} and \eqref{eq:finite_Galois_slc:boundaries}}
}  \\
& = \of^* \left( K_{\oY} + E' + \sum_{\Delta_i \leq E - E'} \left( \frac{1}{r_i} \Delta_i + \frac{r_i -1 }{r_i} \Delta_i \right)  + \sum_{\Delta_i \not\leq E } \frac{r_i -1 }{r_i} \Delta_i \right)\\
& = \of^* \left( K_{\oY} + E  + \sum_{\Delta_i \not\leq E } \frac{r_i -1 }{r_i} \Delta_i \right)\\
\end{split}
\end{equation}
Now, we may conclude the proofs of both points.
\begin{enumerate}
\item  $Y$ is slc if and only if $(\oY,E)$ is log canonical \cite[Defintion-Lemma 5.10]{Kollar_Singularities_of_the_minimal_model_program}. For the latter it is enough to show that $(\oY,E+ \Delta)$ is log canonical for some effective divisor $\Delta$. However, this is shown in \eqref{eq:finite_Galois_slc:computation}, using \cite[Proposition 5.20]{Kollar_Mori_Birational_geometry_of_algebraic_varieties} and that $D' \leq D$.
\item First, note that $\pi^* K_X  \sim K_{\oX} + D$ and $\sigma^* K_Y  \sim K_{\oY} + E$. Second, if $f$ is unramified at the smooth codimension 1 points, then the sum in the last line of \eqref{eq:finite_Galois_slc:computation} is empty. In particular then $\of^* (K_{\oY} + E) = K_{\oX} + D'$. Furthermore, by the other assumption $D= D'$. Hence, the following computation concludes our proof (assuming that we have chosen compatible representatives for $K_X$ and $K_Y$ that avoid the nodal codimension 1 points).
\begin{equation*}
\qquad \quad K_X = \pi_* \pi^* K_X = \pi_*  ( K_{\oX} + D ) =  \pi_* \of^* (K_{\oY} + E) = \pi_* \of^* \sigma^* K_Y = \pi_* \pi^* f^* K_Y = f^* K_Y
\end{equation*}
\end{enumerate}

\end{proof}

\begin{cor}
\label{cor:pullback_K_X}
If $\pi : \sX \to X$ is the coarse moduli morphism introduced in \eqref{eq:stacky_tower} and \eqref{eq:schemey_tower}, then $ \pi^{[*]} \omega_{X/U}^{[N]} \cong \omega_{\sX/U}^{N}$.
\end{cor}

\begin{proof}
It is enough to show the isomorphism in relative codimension one \cite[Proposition 3.6.2]{Hassett_Kovacs_Reflexive_pull_backs}, i.e. we may replace $X$ by its open locus where its fibers  are either smooth or nodal. Then $\omega_{X/U}$ becomes a line bundle. Hence, it is enough to show for the corresponding Cartier divisors that $\pi^* K_{X/U} \sim K_{\sX/U}$. Since $\pi$ is finite and birational, for any compatible choices of $K_{X/U}$ and $K_{Y/U}$, $\pi^* K_{X/U} + E = K_{\sX/U}$  for some divisor $E$ supported on the locus where $\pi$ is not an isomorphism. Furthermore, $E$ is independent of the choice of $K_{X/U}$ and $K_{Y/U}$ as long as they are chosen compatibly. We need to show that $E=0$, or with other words its multiplicity $m_P (E)$ is zero at every codimension one point $P$ of $\sX$. Fix such a $P$. Let $P'$ be the image of $P$ in $X$. For showing that $m_P(E)=0$, we may replace both $\sX$ and $X$ by  \'etale neighborhoods of $P$ and $P'$, respectively. So, by  
\cite[Lemma 2.2.3]{Abramovich_Vistoli_Compactifying_the_space_of_stable_maps} we may assume that $\sX$ is a scheme with an action of a finite 
group $G$, and $X=\sX/G$. Notice first 
that then $\sX$ has demi-normal fibers as well 
by Lemma \ref{lem:stack_demi_normal}. Second, notice that for compatible choice of $K_{\sX/U}$ and $K_{X/U}$, none of the components of $K_{\sX/U} - \pi^* K_{X/U}$ contain a fiber over $U$, because  $\sX \to X$ is isomorphism at the general points of the fibers by Fact \ref{fact:towers}.\ref{itm:towers:general_points}. In particular, it is enough to show that $K_{\sX_u} - \pi_u^* K_{X_u}=0$ for every $u \in U$ and compatible choices of canonical divisors. Third notice that by points \eqref{itm:towers:codimension_one_points_singular} and \eqref{itm:towers:codimension_one_points_regular} of Fact \ref{fact:towers}, Lemma \ref{lem:finite_Galois_slc}.\ref{itm:finite_Galois_slc:unramified} applies and concludes our proof.  
\end{proof}

\begin{lem}
\label{lem:degeneration_Kodaira_fibration}
Consider the commutative diagram
\begin{equation*}
\xymatrix{
X \ar[r]^{\xi} & Y \ar[r]^{\phi} & B
},
\end{equation*}
where
\begin{enumerate}
\item \label{itm:degeneration_Kodaira_fibration:base} $B$ is a smooth curve,
\item \label{itm:degeneration_Kodaira_fibration:first_level} $\phi$ is a family of pure dimensional stacks  of the same dimension at least one,
\item $\xi$ is a flat  family of one dimensional stacks, such that all fibers have only nodal or regular points,
\item \label{itm:degeneration_Kodaira_fibration:second_level_generic} $\xi$ is smooth over $B \setminus \{ 0 \}$,
\item \label{itm:degeneration_Kodaira_fibration:second_level_central} $\xi$ is smooth over a dense set of $Y_0$.
\end{enumerate}
Then, $\xi$ is smooth.
\end{lem}

\begin{proof}
Assume $X$ is not smooth. Choose then a point $P \in X$, such that the fiber $\xi^{-1}(\xi(P))$ is nodal at $P$. By passing to \'etale neighborhoods of $\xi(P)$ and $P$, we may assume that all spaces involved are schemes. However, then we obtain a flat deformation of a node such that the locus of singular fibers has codimension at least two. This is impossible since the node has codimension one in its versal deformation space.
\end{proof}

Finally a lemma used in the proof of Corollary \ref{cor:main_projective}

\begin{lem}
\label{lem:representable}
The map $\pi : \widetilde{\fK\fF}_h \to \fM_h$ is representable in the algebraic space sense. 
\end{lem}

\begin{proof}
By \cite[Lemma 4.4.3]{Abramovich_Vistoli_Compactifying_the_space_of_stable_maps}, the statement of the Lemma is equivalent to saying that for every algebraically closed field $k$ and $x \in \widetilde{\fK\fF}_h(k)$, $\Aut(x) \to \Aut(\pi(x))$ is injective. This follows immediately from the fact that $\pi$ is a forgetful functor, hence every automorphism of $x$ is an automorphism of $\pi(x)$ which respects the extra structure, i.e., the intermediate tower levels, of $x$.
\end{proof}

\section{Arbitrary bases}
\label{sec:arbitrary base}

Here we treat the arbitrary base case. That is, we allow $U$ of Notation \ref{notation:tower} to be affine too. Most of the section is devoted to the proof of Theorem \ref{thm:main_affine}. However,  the proof of Corollary \ref{cor:main_affine} is also included at the end of the section.

First we try to convey an intuition of why considering non-compact bases are much harder then the compact ones. The basic problem is that an entire class of new deformations appear if $U$ is not compact. Intuitively the following happens. Consider a deformation of the tower in Notation \ref{notation:tower}, in the case of $n=2$. Assume for simplicity that the deformation is such that the middle level deforms too. That is, we have a diagram of two Cartesian squares:
\begin{equation*}
\xymatrix{
X=X_2 \ar@/_3pc/[dd]^f \ar[d]^{f_{2}} \ar[r] & X'= X_2' \ar@/^3pc/[dd]^{f'} \ar[d]^{f_{2}'} \\
X_{1} \ar[d]^{f_{1}} \ar[r] & X_{1}' \ar[d]^{f_{1}'} \\
U \ar[r] & U \times T
},
\end{equation*}
where $T$ is a (not necessarily projective) smooth curve. We also assume that $f'$ is smooth. If $U$ was projective, then the smoothness of $f_1$ and $f_2$ and the open property of smoothness would imply that $f_1'$ and $f_2'$ are smooth too. However, as soon as we pass from non-compact to affine base, neither $f_1'$ nor $f_2'$ have any reason to be smooth. In fact, they are not smooth in general. 

Being in a  more subtle situation means, that the proof in this case will be based on a different method. In fact, we prove rigidity using the iterated Kodaira-Spencer map (Definition \ref{def:iterated_Kodaira_Spencer}), as stated in Theorem \ref{thm:main_affine}. In the entire section we work in the situation of Notation \ref{notation:tower}. Since the statement of Theorem \ref{thm:main_affine} is local, we assume that $U$ is affine.  To get rid of the indices, we introduce $Y:=X_1$, $g:=f_2$, $h:=f_1$. Hence we are in the situation
\begin{equation}
\label{eq:affine_case}
\xymatrix{
X \ar[r]^g \ar@/^2pc/[rr]^f & Y  \ar[r]^h & U .
}
\end{equation}
Consider the following  diagram which becomes commutative with exact rows if we disregard the vertical curled arrow.
\begin{equation}
\label{diag:basic}
\xymatrix{
0 \ar[r] & \sT_{X/U} \otimes g^* \sT_{Y/U} \ar[r] & \sT_X \otimes g^* \sT_{Y/U} \ar[r] & f^* \sT_U \otimes g^* \sT_{Y/U} \ar[r] & 0 \\
0 \ar[r] & \sT_{X/U} \otimes \sT_{X/U} \ar[r] \ar[u] \ar[d] & \sT_X \otimes \sT_{X/U} \ar[r] \ar[u] \ar[d] & f^* \sT_U \otimes \sT_{X/U} \ar[r] \ar[u] \ar@{=}[d] & 0 \\
0 \ar[r] & \wedge^2 \sT_{X/U} \ar[r] \ar@/_2pc/[u] & \wedge^2 \sT_X \ar[r] & f^* \sT_U \otimes \sT_{X/U} \ar[r] & 0
},
\end{equation}
where
\begin{itemize}
\item the homomorphism $\sT_{X/U} \otimes \sT_{X/U} \to \wedge^2 \sT_{X/U}$  is the  wedge product map,
\item the homomorphism $\sT_{X} \otimes \sT_{X/U} \to \wedge^2 \sT_X$ is the embedding $\sT_{X} \otimes \sT_{X/U} \to \sT_{X} \otimes \sT_{X}$ composed with the wedge product $\sT_{X} \otimes \sT_{X} \to \wedge^2 \sT_{X}$ and
\item the homomorphism $\wedge^2 \sT_{X/U} \to  \sT_{X/U} \otimes \sT_{X/U}$ is the splitting of the wedge product map, given by $a \wedge b \mapsto \frac{1}{2} ( a \otimes b - b \otimes a)$
\end{itemize}
Recall the homomorphisms $\rho_i$ from Definition \ref{def:iterated_Kodaira_Spencer}. 
Our aim is to show that the map   $\rho_2 : R^1 f_* (f^* \sT_U \otimes \sT_{X/U}) \to R^2 f_* (\wedge^2 \sT_{X/U})$ is injective on the image of   $\rho_1 : \sT_U^{\otimes 2} \to R^1 f_* (f^* \sT_U \otimes \sT_{X/U})$. Clearly, that will yield the injectivity of $\ks_f = \rho_2 \circ \rho_1$.
\begin{notation}
\label{notation:main_diagram}
Taking long exact sequences of derived pushforwards of the rows of  \eqref{diag:basic} yields the following diagram. We also introduce names for certain homomorphisms in the diagram.
\begin{equation*}
\xymatrix{ 
R^1 f_* ( g^* \sT_{Y/U} \otimes f^* \sT_U) \ar[r]^{\eta} & R^2 f_* (\sT_{X/U} \otimes g^* \sT_{Y/U}) \\
R^1 f_* (\sT_{X/U} \otimes f^*\sT_U) \ar[r]^{\beta} \ar@{=}[d] \ar[u]^{\delta} & R^2 f_* ( \sT_{X/U} \otimes \sT_{X/U}) \ar[u] \ar[d]^{\gamma} \\
R^1 f_* (\sT_{X/U} \otimes f^*\sT_U) \ar[r]^{\alpha} & R^2 f_* (  \wedge^2 \sT_{X/U}) \ar@/_2pc/[u]^{\varepsilon}
}
\end{equation*}
\end{notation}

Now we prove Theorem \ref{thm:main_affine}. In fact, important parts of the proof are done in Propositions \ref{prop:delta}, \ref{prop:eta} and \ref{prop:epsilon}, afterwards.

\begin{proof}[Proof of Theorem \ref{thm:main_affine}] 
We use Notation \ref{notation:main_diagram}. By Proposition \ref{prop:delta} and Proposition \ref{prop:eta} both $\delta$ and $\eta$ are generically injective. Hence, so is $\beta$. Consider now the following diagram.
\begin{equation*}
\xymatrix{
 \sT_U^{\otimes 2} \ar[r]_(0.3){\ks_f=\rho_1} \ar@{=}[d] \ar@/^2pc/[rr]^{\nu} &  \sT_U \otimes R^1 f_*( \sT_{X/U}) \ar[r]^{\beta} \ar@{=}[d] & R^2 f_*(  \sT_{X/U} \otimes \sT_{X/U} ) \ar[d]^{\gamma} \\
 \sT_U^{\otimes 2} \ar[r]^(0.3){\ks_f=\rho_1} \ar@/_2pc/[rr]_{\iks_f} &  \sT_U \otimes R^1 f_*( \sT_{X/U}) \ar[r] & R^2 f_*( \wedge^2  \sT_{X/U}  ) \ar@/_2pc/[u]^{\varepsilon}
} 
\end{equation*}
Since $h$ has variation 1, the same holds for $f$. One reason is for example that a variety has only finitely many dominant general type images up to birational equivalence (e.g., 
\cite[Corollary 1.4]{Hacon_McKernan_Boundedness_of_pluricanonical_maps_of_varieties_of_general_type}).  Hence $\ks_f$ is injective. Then  $\nu := \beta \circ \ks_f$ is generically injective and also injective, since for homomorphisms from torsion free sheaves on varieties generic injectivity implies injectivity. By Proposition \ref{prop:epsilon}, $\im \nu \subseteq \im \varepsilon$. Since $\varepsilon$ is a splitting of the surjection $\gamma$, this means, that $\gamma$ maps $\im \nu$ injectively. Hence $\iks_f = \gamma \circ \nu$ is injective too.
\end{proof}

The rest of the section deals with the propositions referenced by the proof of Theorem \ref{thm:main_affine}.

\begin{prop}
\label{prop:delta}
In the situation of Notation \ref{notation:main_diagram}, $\delta$ is generically injective.
\end{prop}

\begin{proof}
Consider the following exact sequence.
\begin{equation}
\label{eq:tangent}
\xymatrix{ 
0 \ar[r] & \sT_{X/Y} \ar[r] & \sT_{X/U} \ar[r] & g^* \sT_{Y/U} \ar[r] & 0 
}
\end{equation}
Since $g$ has maximal variation, $X_u \to Y_u$ is non-isotrivial for generic $u \in U$. Hence for generic $u \in U$, $\omega_{X_u/Y_u}$ is ample by Corollary \ref{cor:rel_canonical_sheaf_big}. Then by Kodaira vanishing
$H^1(X_u,\sT_{X/Y})=0$. So, $R_1 f_* \sT_{X/Y}$ is torsion. Hence, taking the long exact sequence of derived pushforwards of \eqref{eq:tangent} yields that  the natural map 
\begin{equation*}
R^1f_* \sT_{X/U} \to  R^1 f_* g^* \sT_{Y/U}
\end{equation*} 
is generically an injection.
\end{proof}

For the next proposition we need a lemma first.

\begin{lem}
\label{lem:ample}
If $\sE$ is an ample vector bundle over a projective smooth curve, $\xi : \sE \to \sH$ a generically surjective homomorphism onto a vector bundle, then $\ker \xi \otimes \det \sH$ is an ample vector bundle.
\end{lem}

\begin{proof}
First, since $\im \xi  \subseteq \sH$ is of full rank, $\det (\im \xi)$ is a sub-line bundle of $\det \sH$. Hence, it is enough to show that $\ker \xi \otimes \det ( \im \xi)$ is ample. With other words, by replacing $\sH$ with $\im \xi$, we may assume that $\xi$ is surjective. 

Assume now that there is a surjection $\phi: \ker \xi \otimes \det \sH \to \sM$ onto a line bundle. Then one can form the pushout diagram
\begin{equation}
\label{eq:ample:pushforward}
\xymatrix{
0 \ar[r] & \ker \xi \otimes \det \sH \ar[r] \ar[d]^{\phi} & \sE \otimes \det \sH \ar[r] \ar[d] & \sH \otimes \det \sH  \ar@{=}[d]  \\
0 \ar[r] & \sM \ar[d] \ar[r] & \sF \ar[d] \ar[r] & \sH \otimes \det \sH  
\\
& 0  & 0 ,
}
\end{equation}
where $\sF:= \factor{\sE \otimes \det \sH }{\ker \phi}$. Since $\sE$ is ample, so is $\sF \otimes (\det \sH)^{-1}$. That is, 
\begin{multline*}
\deg \sF > (\rk \sF) \deg (\det \sH) = (\rk \sH +1) \deg ( \det \sH) \\ = \det ( \sH \otimes \det \sH).  
\end{multline*}
This implies, by the bottom exact row of \eqref{eq:ample:pushforward}  that $\deg \sM > 0$ (notice, that by construction the right most edge in that row is surjective). Hence all line bundle quotients of $\ker \xi \otimes \det \sH$ have positive degree. If $\tau$ is a  finite map of smooth curves, then the same holds for $ \tau^* (\ker \xi \otimes \det \sH)$, since it is isomorphic to $(\ker (\tau^* \sE \to \tau^* \sH)) \otimes \tau^* \det \sH$ and $\tau^* \sE$ is ample too. This shows that $\ker \xi \otimes \det \sH$ is indeed ample.
\end{proof}

\begin{prop}
\label{prop:eta}
In the situation of Notation \ref{notation:main_diagram}, $\eta$ is generically injective. 
\end{prop}

\begin{proof}
To prove the generic injectivity of $\eta$ we would need that
\begin{equation*}
R^1 f_* (\sT_X \otimes g^* \sT_{Y/U} ) 
\end{equation*}
is torsion. First, we show that 
\begin{equation}
\label{eq:g_push_forward}
g_* ( \sT_X \otimes g^* \sT_{Y/U}) \cong  g_* \sT_X \otimes \sT_{Y/U}  = 0  .
\end{equation}
Consider the following exact sequence.
\begin{equation*}
\xymatrix{
0 \ar[r] & \sT_{X/Y} \ar[r] & \sT_X \ar[r]  & g^* \sT_Y \ar[r] & 0
}
\end{equation*}
Then by the pushfoward long exact sequence we obtain
\begin{equation*}
\xymatrix{
0 \ar[r] & g_* \sT_{X/Y} \otimes \sT_{Y/U} = 0 \ar[r] & g_*\sT_X \otimes \sT_{Y/U} & \\  & \ar[r] & \sT_Y \otimes \sT_{Y/U} \ar[r]   &  R^1 g_* \sT_{X/Y} \otimes \sT_{Y/U} ,
}
\end{equation*}
where the last map is $\ks_g$ tensored with  $\sT_{Y/U}$.
Since $g$ has maximal variation, this map is injective, which proves (\ref{eq:g_push_forward}).

So, by the Grothendieck spectral sequence it is enough to show that
\begin{equation*}
h_* R^1 g_* (\sT_X \otimes g^* \sT_{Y/U} )
\end{equation*}
is torsion. By relative duality the following isomorphisms hold.
\begin{equation*}
h_* R^1 g_*(\sT_{X} \otimes g^* \sT_{Y/U}) \cong  h_* (( g_*( \omega_{X/Y} \otimes \Omega_{X} \otimes g^* \omega_{Y/U}) )^*)  \cong h _* (( g_*( \omega_{X/U} \otimes \Omega_{X} ) )^*)
\end{equation*}
So it is enough to show, that $h^0(Y_u, (g_*( \omega_{X/U} \otimes \Omega_{X} ))^*)=0$ for generic $u \in U$. In general $g_*( \omega_{X/U} \otimes \Omega_{X} )$ is not a vector bundle, however being the pushforward of a torsion free sheaf it is torsion free at least. Hence, since $Y$ is smooth and $\dim Y=2$, it is a vector bundle except at finitely many points. By leaving out the images of those points from $U$, we may assume, that $g_*( \omega_{X/U} \otimes \Omega_{X} )$ is in fact locally free.

Consider the following short exact sequence.
\begin{equation*}
\xymatrix{
0 \ar[r] & \omega_{X/U} \otimes g^* \Omega_{Y} \ar[r] & \omega_{X/U} \otimes \Omega_{X} \ar[r] & \omega_{X/U} \otimes \omega_{X/Y} \ar[r] & 0
}
\end{equation*}
By pushing it forward one obtains
\begin{equation*}
\xymatrix{
0 \ar[r] & g_* \omega_{X/Y} \otimes \omega_{Y/U} \otimes \Omega_{Y} \ar[r] & g_* (\omega_{X/U} \otimes \Omega_{X}) & \\
  \ar[r] & g_* (\omega_{X/Y}^{\otimes 2})   \otimes \omega_{Y/U} \ar[r]^(0.4){\mu} &   R^1 g_* \omega_{X/U} \otimes \Omega_{Y} \cong \omega_{Y/U} \otimes \Omega_Y ,
}
\end{equation*}
where the last homomorphism $\mu$ is the dual of $\ks_g$  tensored with $\omega_{Y/U}$. Again, $\im \mu$ is  not necessarily locally free, but being a subsheaf of a torsion free sheaf, it is torsion free. Hence as before, we may assume that it is in fact locally free. Then we see, that $g_*( \omega_{X/U} \otimes \Omega_{X} )$ is the extension of two locally free sheaves: $ g_* \omega_{X/Y} \otimes \omega_{Y/U} \otimes \Omega_{Y}$ and  $\ker \mu$. We conclude our proof, by showing that 
\begin{equation}
\label{eq:vanishing_first}
 h^0(Y_u, (g_* \omega_{X/Y} \otimes \omega_{Y/U} \otimes \Omega_{Y})^*) = 0
\end{equation}
and 
\begin{equation}
\label{eq:vanishing_second}
 h^0(Y_u, (\ker \mu)^*) = 0
\end{equation}
for generic $u \in U$. 

For the first one, since $\Omega_Y$ is the extension of two nef line bundles by Proposition \ref{prop:rel_canonical_nef}, it is also nef.  The pushforward $g_* \omega_{X/Y}$ is nef too  (e.g., \cite[Theorem 4.1]{Viehweg_Weak_positivity}) and $\omega_{Y/U}$ is $h$-ample. Hence $g_* \omega_{X/Y} \otimes \omega_{Y/U} \otimes \Omega_{Y}$ is ample on $Y_u$ for each $u \in U$. This implies  \eqref{eq:vanishing_first} for every $u \in U$.

To show \eqref{eq:vanishing_second}, notice that $\mu$ is generically surjective by the assumption $\var g =2$. By possibly restricting $U$ we may assume that $\mu$ is generically surjective on each $Y_u$. Notice, that  $h^0(X_y, \omega_{X_y}^{\otimes 2})$ is a constant function of $y$. Hence by \cite[Corollary III.12.9]{Hartshorne_Algebraic_geometry},  $g_* (\omega_{X/Y}^{\otimes 2})$ is locally free. Fix some $u \in U$. If we restrict $\mu$ to $Y_u$, using that $\im \mu$, $\ker \mu$ and $f_* (\omega_{X/Y}^{\otimes 2})$ are locally free, we obtain the exact sequence
\begin{equation*}
\xymatrix{
0 \ar[r] &  (\ker \mu)|_{Y_u} \ar[r] & g_* (\omega_{X/Y}^{\otimes 2})   \otimes \omega_{Y_u} \ar[r]^(0.6){\mu|_{Y_u}}  &  \Omega_Y \otimes \omega_{Y_u}  
},
\end{equation*}
where the last map is generically surjective. Then, since $\det (\Omega_Y|_{Y_u}) \cong \omega_{Y_u}$, by Lemma \ref{lem:ample}, Proposition \ref{prop:rel_canonical_nef} and Lemma \ref{lem:positivity_pushforwards}, $(\ker \mu)|_{Y_u}$ is ample. Then \eqref{eq:vanishing_second} follows, since by $\ker \mu$ being locally free, $(\ker \mu)^*|_{Y_u} \cong ((\ker \mu)|_{Y_u})^*$.
\end{proof}

\begin{prop}
\label{prop:epsilon}
The image of the composition 
\begin{equation*}
\xymatrix{
 \sT_U^{\otimes 2} \ar[r] \ar@/^2pc/[rr]^{\nu} &  \sT_U \otimes R^1 f_*( \sT_{X/U}) \ar[r] & R^2 f_*(  \sT_{X/U} \otimes \sT_{X/U} )
}
\end{equation*}
is contained in $\im \varepsilon$ (see Notation \ref{notation:main_diagram} for the definition of $\varepsilon$).
\end{prop}

\begin{proof}
Since $\sT_U^{\otimes 2}$ is a line bundle and $U$ is affine, by possibly restricting $U$ we may assume that $\sT_U \cong \sO_U$. This yields a generator $t$ of $\sT_U$. We have to show that $\nu( t \otimes t) \subseteq \im \varepsilon$. 

Since we assumed that $U$ is affine, we can replace the derived pushforwards by global cohomology in Notation \ref{notation:main_diagram}. Then we get the following diagram.
\begin{equation*}
\xymatrix{
 H^0(X, f^* \sT_U \otimes f^* \sT_U) \ar[r]^{\ks_f} \ar[dr]^{\nu} & H^1(X,  f^* \sT_U \otimes \sT_{X/U} ) \ar[d]^{\beta} \\ & H^2(X,  \sT_{X/U} \otimes \sT_{X/U} )
}
\end{equation*}
We are going to use Dolbeault cohomology to prove the statement of the proposition. Let $\bar{t}$ be the element of $H^0(X, f^* \sT_U)$ corresponding to $t \in \sT_U$. Then there is an element $a \in \sA^{0,0}(f^* \sT_U)$ in the Dolbeault resolution corresponding to $\bar{t}$. Let $b \in  \sA^{0,0}(\sT_X)$ be any lift of $a$, and $c := \bar{\partial} b$. By abuse of notation we will view $c$ both as an element of $\sA^{0,1}(\sT_X)$ and of $\sA^{0,1}(\sT_{X/U})$. Because of the presence of two different wedge products (one on antiholomorphic forms, and one one tangent bundles), we will need to write $c$ in local coordinates:
\begin{equation}
\label{eq:local}
 c := \sum_{i=1}^3 c_i d\bar{z}_i  \qquad c_i \in \sT_{X/U} 
\end{equation}
First we compute $\nu(a \otimes a)$. To get $\ks_f(a \otimes a)$, we have to compute a boundary homomorphism of the exact sequence
\begin{equation*}
\xymatrix{
0 \ar[r] & f^* \sT_U \otimes \sT_{X/U}  \ar[r] &  f^* \sT_U \otimes \sT_X \ar[r] & f^* \sT_U \otimes f^* \sT_U \ar[r] & 0
}
\end{equation*}
So, we lift $a \otimes a$ to get $a \otimes b$ and then we apply $\bar{\partial}$. We obtain the following (remember, $a$ is holomorphic, hence $\bar{\partial}(a)=0$).
\begin{equation*}
\nu(a \otimes a) = \bar{\partial}( a \otimes b) = a \otimes \bar{\partial}(b) = a \otimes c \in \sA^{0,1}(\sT_{X/U} \otimes f^* \sT_U)
\end{equation*}
That is,
\begin{equation*}
\nu(a \otimes a)=\beta (\ks_f(a \otimes a)) = \beta ( c \otimes a) 
\end{equation*}
is obtained by feeding $a \otimes c$ to the edge morphism of the exact sequence
\begin{equation*}
\xymatrix{
0 \ar[r] &  \sT_{X/U} \otimes \sT_{X/U} \ar[r] & \sT_X \otimes \sT_{X/U} \ar[r] & f^* \sT_U \otimes \sT_{X/U} \ar[r] & 0 .
}
\end{equation*}
I.e., we lift $a \otimes c$ to $b \otimes c$, and then we apply $\bar{\partial}$. We obtain the following (remember, $c$ is Dolbeault-closed, hence $\bar{\partial}(c)=0$).
\begin{equation*}
 \nu(a \otimes a) = \bar{\partial}( b \otimes c) = \bar{\partial}(b) \otimes c = c \otimes c \in \sA^{0,2}(\sT_{X/U} \otimes  \sT_{X/U})
\end{equation*}
 We conclude our proof by showing that $\varepsilon (\gamma(c \otimes c))=c \otimes c$. This part is slightly confusing, so we change to the local expression of \eqref{eq:local}. 
\begin{eqnarray*}
\varepsilon (\gamma(c \otimes c))&  = & \varepsilon \left(\gamma \left(\sum (c_i \otimes c_j) d\bar{z}_i \wedge d \bar{z}_j \right) \right)
\\& = & \varepsilon \left( \sum (c_i \wedge c_j) d\bar{z}_i \wedge d \bar{z}_j \right)
\\ & = & \sum \frac{1}{2}\left( c_i \otimes c_j - c_j \otimes c_i \right) d\bar{z}_i \wedge d \bar{z}_j
\\ & = &  \frac{1}{2} \sum_{i < j} ( (  c_i \otimes c_j - c_j \otimes c_i ) d\bar{z}_i \wedge d \bar{z}_j \\ & & \qquad + ( c_j \otimes c_i - c_i \otimes c_j ) d\bar{z}_j \wedge d \bar{z}_i )
\\ & = &   \sum_{i < j}  ( c_i \otimes c_j - c_j \otimes c_i ) d\bar{z}_i \wedge d \bar{z}_j 
\\ & =  &  \sum   c_i \otimes c_j  d\bar{z}_i \wedge d \bar{z}_j  = c \otimes c
\end{eqnarray*}

\end{proof}

\begin{proof}[Proof of Corollary \ref{cor:main_affine}]
Assume we are in the situation described in Corollary \ref{cor:main_affine}. Then $\iks_{\nu \circ \tau} $ is injective. We want to prove that $\iks_{\nu}$ is injective as well. Let $f : X \to U$ and $f' : X' \to U'$ be the families associated to $\nu$ and $\nu \circ \tau$, respectively. In particular, then $X' \cong X \times_U U'$ and $\sT_{X'/U'} \cong (\tau')^* \sT_{X/U} $, where $\tau'$ is the natural projection map $X' \to X$. 

For any $1 \leq p \leq n$, let 
\begin{equation*}
0=\sF_0^p \subseteq \sF_1^p \subseteq \dots \subseteq \sF_{p}^p \subseteq \sF_{p+1}^p=\wedge^p \sT_X 
\end{equation*}
and 
\begin{equation*}
0=\widetilde{\sF}_0^p \subseteq \widetilde{\sF}^p_1 \subseteq \dots \subseteq \widetilde{\sF}^p_{p} \subseteq \widetilde{\sF}^p_{p+1}=\wedge^p \sT_{X'} 
\end{equation*}
be the filtrations of Definition \ref{def:iterated_Kodaira_Spencer} for $g=f$ and $g=f'$, respectively. From the construction of the filtrations $\sF_{\bullet}^p$ and $\widetilde{\sF}^p_{\bullet}$, one can see that there are natural morphisms $\phi_i : \widetilde{\sF}_i^p \to \tau^* \sF^p_i$, such that  $\phi_i|_{\widetilde{\sF}_j^p} = \phi_j$ for every $j< i$. Furthermore, $\phi_i$ induce the natural homomorphisms 
\begin{equation*}
 ((f')^* \wedge^i \sT_{U'}) \otimes  (\wedge^{p-i} \sT_{X'/U'}) \to  ( (\tau')^* f^* \wedge^i  \sT_U) \otimes  ( (\tau')^* \wedge^{p-i} \sT_{X/U})
\end{equation*}
on the quotients via the isomorphism \eqref{eq:iterated_Kodaira_Spencer:quotient}. Then, the short exact sequence \eqref{eq:iterated_Kodaira_Spencer:short_exact} is compatible with $\tau$ in the sense that there is a commutative diagram as follows.
\begin{equation}
\label{eq:main_affine:compatibility}
\xymatrix{
 \wedge^p \sT_{X'/U'}  \ar[r] \ar[d]
& \widetilde{\sF}^p_2   \ar[r] \ar[d]
& (f')^* \sT_{U'} \otimes \wedge^{p-1} \sT_{X'/U'}   \ar[d] 
\\
(\tau')^* \wedge^p  \sT_{X/U}  \ar[r] 
& (\tau')^* \sF_2 \ar[r]   
& ((\tau')^* f^* \sT_U) \otimes ((\tau')^* \wedge^{p-1} \sT_{X/U} ) 
}
\end{equation}

At this point, we have to make a short digression. Let $\sG$ be any of the finitely many sheaves on $X$ that has appeared  in the proof. Generically on $U$,  $R_i f_* (\sG)$ commutes with base change by \cite[Theorem III.12. 8 and Corollary III.12.9]{Hartshorne_Algebraic_geometry}. Hence, by restricting $U'$ we might as well assume that $R^i f_*( \sG)$ commutes with base change via $\tau $ for all the mentioned sheaves $\sG$. In particular, taking the long exact sequence of derived pushforwards of the bottom line in \eqref{eq:main_affine:compatibility} is compatible with base change via $\tau$. I.e.,  \eqref{eq:main_affine:compatibility} yields the following compatibilities of edge morphisms considered in Definition \ref{defn:iterated_Kodaira_fibration}.
\begin{equation*}
\xymatrix{
\sT_{U'}^{\otimes (n-p+1)} \otimes R^{p-1} (f')_*( \wedge^{p-1} \sT_{X'/U'}) \ar[r]^(0.55){\rho_p'} \ar[d] 
& \sT_{U'}^{\otimes (n-p)} \otimes R^{p} (f')_*  (\wedge^p \sT_{X'/U'}) \ar[d] \\
\tau^* ( \sT_{U}^{\otimes (n-p+1)} \otimes R^{p-1} f_*(  \wedge^{p-1} \sT_{X/U}) ) \ar[r]^(0.55){\tau^* \rho_p} 
& \tau^* ( \sT_{U}^{\otimes (n-p)} \otimes R^{p} f_*  ( \wedge^p \sT_{X/U}) )
}
\end{equation*}
However, then we obtain similar compatibility for the composition of these edge maps as follows.
\begin{equation*}
\xymatrix{
\sT_{U'}^{\otimes n} \ar[r]^(0.3){\iks_{\nu \circ \tau}} \ar[d] &  R^{n} (f')_*  (\wedge^n \sT_{X'/U'}) \ar[d]^{\cong} \\
\tau^* \sT_{U}^{\otimes n}  \ar[r]^(0.35){\tau^* \iks_{\nu}} &  \tau^* R^{n} f_*  (\wedge^p \sT_{X/U}) 
}
\end{equation*}
Notice first that the right vertical arrow is an isomorphism. Second, the homomorphism $\sT_{U'}^{\otimes n} \to \tau^* \sT_{U}^{\otimes n}$ is isomorphism, since $\tau$ is \'etale. Hence, $\iks_{\nu}$ is injective. Then, rigidity of $\nu$ follows via Theorem \ref{thm:iks_injective_rigid}.
\end{proof}

\section{Rigidity and variational rigidity}
\label{sec:variationally_rigid}

This is a minor addendum about how \cite[Corollary 8.4]{Viehweg_Zuo_Discreteness_of_minimal_models_of_Kodaira_dimension_zero_and_subvarieties_of_moduli_stacks} and \cite[Theorem 4.14]{Kovacs_Strong_non_isotriviality_and_rigidity} implies Theorem \ref{thm:iks_injective_rigid}. In the mentioned papers, the following rigidity property is shown for morphisms with injective iterated Kodaira-Spencer morphisms. 

\begin{defn}
\label{def:variationally_rigid}
A morphism $\nu : U \to \fM_h$ from a manifold is \emph{variationally rigid}, if for every deformations $\nu' : U \times T \to \fM_h$ over a smooth irreducible curve $T$, $\dim (\im \nu') = \dim (\im \nu)$.
\end{defn}

However, this rigidity property a priori is slightly weaker than of Definition \ref{def:rigid}. Luckily, in reality, at least for one dimensional bases, the two rigidity properties are equivalent as stated in Proposition \ref{prop:variationally_rigid}. 

Before starting the argument, note that by usual base extension and restriction arguments we may assume that we work over the complex numbers. Then we may also use the holomorphic category at certain points. However, there will be some indication whenever we do that. With other words,  every space and map is algebraic unless otherwise stated.

\begin{prop}
\label{prop:variationally_rigid}
If a finite morphism $\nu : U \to \fM_h$ from a smooth curve is variationally rigid, then it is rigid.
\end{prop}

The main ingredient in the proof of Proposition \ref{prop:variationally_rigid} is the uniformization of one dimensional smooth Deligne-Mumford stacks, i.e. of smooth DM-curves, worked out in \cite{Behrend_Noohi_Uniformization_of_Deligne_Mumford_curves}. There, it is  shown that the possible universal covers of smooth DM-curves are the hyperbolic plane $\bH$, the complex plane $\bC$ or the weighted projective line $\bP(m,n)$ for  $(m,n)=1$.
First, we show the Shafarevich conjecture over DM-curve bases in Lemma \ref{lem:stacky_shafarevich}. 

\begin{defn} \cite[ Definition 7.3, Corollary 7.7]{Behrend_Noohi_Uniformization_of_Deligne_Mumford_curves}
\label{def:hyperbolic}
A smooth \emph{Deligne-Mumford curve}, is \emph{hyperbolic} if its universal cover is $\bH$. Equivalently, for algebraic DM-curve $U$, if there is a smooth algebraic curve $V$ with an action of a finite group $G$ such that $U=[V/G]$ and for any such $V$, if $C$ is the smooth compactification of $V$, then $\deg \omega_C((C \setminus V)_{\red}) >0$. 
\end{defn}

\begin{lem}
\label{lem:stacky_shafarevich}
If $\nu : U \to \fM_h$ is a finite morphism from a  smooth Deligne-Mumford curve, then $U$ is hyperbolic. 
\end{lem}

\begin{proof}
Assume $U$ is not hyperbolic. Then by \cite[Proposition 7.2]{Behrend_Noohi_Uniformization_of_Deligne_Mumford_curves}, its universal cover $\widetilde{U}$ contains an open set $\widetilde{U}^0$ which is isomorphic to $\bC \setminus \{0\}$.  Hence, there is a holomorphic map $\bC \setminus \{0\} \to \fM_h$ with one dimensional image. Composing this with the exponential map one obtains a holomorphic map $\bC\to \fM_h$ with one dimensional image. However, this is impossible by \cite[Theorem 0.1]{Viehweg_Zuo_On_the_Brody_hyperbolicity_of_moduli_spaces_for_canonically_polarized_manifolds}.
\end{proof}

Second, we need a stacky version of the De Franchis Theorem shown in Lemma \ref{lem:rigidity_maps_hyperbolic_stack_curves}. The proof of Lemma \ref{lem:rigidity_maps_hyperbolic_stack_curves} also uses the algebraicity of certain holomorphic maps and some properties of torsors  shown in Lemmas \ref{lem:holomorphic_algebraic} and \ref{lem:etale_torsor_trivialization}. 

\begin{lem}
\label{lem:holomorphic_algebraic}
If $X$ and $Y$ are smooth projective curves, $U \subseteq X$ is a (Zariski) open set and $f : U \to Y$ a holomorphic map with finite fibers, then $f$ is algebraic.
\end{lem}

\begin{proof}
By the GAGA principle, it is enough to show that $f$ extends to $X$. For that choose a point $P \in Y \setminus f(U)$. If a priori there is no such point, then by discarding finitely many points of $U$, we may find one.  Choose an (algebraic) embedding $\iota : Y \hookrightarrow \bP^N$, such that $Y$ is not contained in any hyperplane, and there is a hyperplane $H \subseteq \bP^n$, for which $(H \cap Y)_{\red}= P$. Let $x^1, \dots, x^N$ be the coordinate functions of $\bP^n \setminus H$. Then, $x^i \circ f$ are holomorphic. Furthermore, by the finiteness assumption of the lemma, and the non-degeneracy of $\iota(Y)$, $x^i \circ f$ has finite fibers. Hence, every $x^i \circ f$ has non-essential singularity at every point of $X \setminus U$. However, then $[1,x^1 \circ f, \dots , x^N \circ f]$ gives an extension of $\iota \circ f$ over $X$, and consequently of $f$ as well.
\end{proof}

\begin{cor}
\label{cor:holomorphic_algebraic}
If $X$ and $Y$ are projective curves, $U \subseteq X$ and $V \subseteq Y$ (Zariski) open subsets and $f : U \to V$ is a finite level covering map (or in particular a biholomorphism), then $f$ is algebraic.
\end{cor}

\begin{lem}
\label{lem:etale_torsor_trivialization}
If $P \to V \times T$ is an $H$ torsor, where $V$ and $T$ are smooth curves and $H$ is a finite group, then there is a finite \'etale cover $\widetilde{T} \to T$, such that for any $t \in T$, $P \times_T \widetilde{T} \cong P_t \times \widetilde{T}$ as $H$-torsors over $V \times \widetilde{T}$. 
\end{lem}

\begin{proof}
First, we claim that for every $t, t' \in T$, $P_t \cong P_{t'}$ as schemes over $V$. By Corollary \ref{cor:holomorphic_algebraic},  it is enough to show this isomorphism as holomorphic coverings, and then as topological coverings of $V$. However, then to prove the claim, we may replace $T$ by the unit disk $D$, in which case the strong deformation equivalence of $V$ and $V \times T$ concludes the claim.
Similarly, we may prove that for every $v, v' \in T$, $P_v \cong P_{v'}$ as schemes over $T$.

Define then $\widetilde{T}:=P_v$ for arbitrary $v \in V$, $\widetilde{V}:=P_t$ for arbitrary $t \in T$ and $P':=P \times_T \widetilde{T}$. Consider $\Isom_{V \times \widetilde{T}}(P', \widetilde{V} \times \widetilde{T})$. It inherits an $H$-torsor structure from $\widetilde{V} \times \widetilde{T}$. For any $v \in V$, 
\begin{equation*}
 P'|_{\{v\} \times \widetilde{T}} \cong \widetilde{T} \times_T \widetilde{T} 
\underbrace{\cong H \times \widetilde{T}}_{\parbox{130pt}{\tiny the diagonal map $\widetilde{T} \to \widetilde{T} \times_T \widetilde{T}$ is a section, hence $\widetilde{T} \times_T \widetilde{T}$ is a trivial torsor}}
\cong \widetilde{V} \times \widetilde{T}|_{\{v\} \times \widetilde{T}},
\end{equation*}
and for any $ t \in \widetilde{T}$,
\begin{equation*}
 P'|_{V \times \{t\}} \cong 
\underbrace{P|_{V \times \{t'\}}}_{\textrm{$t'$ is the image of $t$ in $T$}} 
\cong \widetilde{V}
\cong \widetilde{V} \times \widetilde{T}|_{V \times \{t\}}.
\end{equation*}
Hence the restriction of $\Isom_{V \times \widetilde{T}}(P', \widetilde{V} \times \widetilde{T})$ to any slice, in any of the two directions, is a trivial $H$-torsor. That is, the action of the normal subgroups $
\pi_1(V,v)$ and $\pi_1(\widetilde{T},t) \subseteq \pi_1(V \times \widetilde{T},(v,t))  $
on $\Isom_{V \times \widetilde{T}}(P', \widetilde{V} \times \widetilde{T})$  is trivial, for arbitrary choice of basepoints. Then so is the action of $\pi_1(V \times \widetilde{T},(v,t))\cong \pi_1(V,v) \times \pi_1(\widetilde{T},t)$. Hence, $\Isom_{V \times \widetilde{T}}(P', \widetilde{V} \times \widetilde{T})$ is the trivial $H$-torsor, that is, there is a global isomorphism as stated in the lemma. This concludes our proof.
\end{proof}

Recall that a morphism $f : U \to W$ is \emph{rigid} if for every deformation $f : U \times T \to W$ over a smooth, irreducible curve $T$, $f_t = f$ for every $t \in T$.
\begin{lem}
\label{lem:rigidity_maps_hyperbolic_stack_curves}
If $f : U \to W$ is a finite morphism between hyperbolic smooth  Deligne Mumford curves, then $f$ is rigid.
\end{lem}

\begin{proof}
Assume $f$ is not rigid. I.e. it has a deformation $f' : U \times T \to W$, where $T$ is a smooth, irreducible curve, and $f'_t \neq f$ for some $t \in T$. By applying Definition \ref{def:hyperbolic}, there is a smooth curve $V$ and a group action of a finite group $G$ on $V$, such that $U= [V/G]$. This yields a map $g : V \to W$ and its deformation $g' : V \times T \to W$. Applying Definition \ref{def:hyperbolic} to $W$ yields a smooth curve $Z$ with a finite group $H$ acting on it, such that $W= [Z/H]$. By the definition of stack quotient then the maps $f$ and $f'$ correspond to $H$-torsors $P$ and $P'$ over $V$ and $V \times T$, respectively, endowed with maps $h$ and $h'$ to $Z$. The situation is summarized in  the following diagram. 
\begin{equation*}
\xymatrix{
&  &  P' \ar[ddrrr]^{h'} \ar@{.>}[d] \\
& P \ar@{.>}[d] \ar[drrrr]^(0.7)h \ar@{-->}[ur] &  V \times T \ar[ddrrr]_(0.4){g'} \ar@{.>}[d] \\
& V \ar@{.>}[d] \ar[drrrr]_(0.4)g \ar@{-->}[ur] &  U \times T \ar[drrr]^(0.4){f'} & & & Z \ar@{.>}[d] \\
[V/G] \ar@{=}[r]  & U \ar[rrrr]^(0.3)f \ar@{-->}[ur] & & & & W \ar@{=}[r] & [Z/H]
}
\end{equation*}
By Lemma \ref{lem:etale_torsor_trivialization}, there is an \'etale cover $\xi : \widetilde{T} \to T$, such that $P' \times_{T} \widetilde{T} \cong P \times \widetilde{T}$ as schemes over $V \times \widetilde{T}$.
Let $h''$ be the map $P' \times_T \widetilde{T} \to Z$ given by the composition of the projection $P' \times_T \widetilde{T} \to P'$  with $h'$. Then  $h''$ is a deformation of $h$ over $\widetilde{t}$. Hence $(h'')_t=h$ for all $t \in \widetilde{T}$. However, by our starting assumption $f'_t \neq f $ for some $t \in T$, which implies that $h_t'' \neq h$ for some $t \in T$. This is a contradiction. Hence our assumption was false, and $f$ is rigid.
\end{proof}

\begin{proof}[Proof of Proposition \ref{prop:variationally_rigid}]
By Lemma \ref{lem:stacky_shafarevich}, $U$ is a hyperbolic curve. Define $W$ to be the normalization of $\im \nu$. Then $W$ is a smooth DM-curve. Hence by Lemma \ref{lem:stacky_shafarevich}, it is also hyperbolic. 

Take any deformation $\nu': U \times T \to \fM_h$ of $\nu$ over a smooth, irreducible curve $T$. Since $\nu$ is variationally rigid, $\nu'$ factors through $W$. That is, there is a commutative diagram as follows.
\begin{equation*}
\xymatrix{
U \times T \ar@/^2pc/[rr]^{\nu'} \ar[r]^{\xi'} & W \ar[r]^{\iota} & \fM_h \\
U \ar[u] \ar@/_2pc/[rr]_{\nu} \ar[r]^{\xi} & W \ar@{=}[u] \ar[r]_{\iota} & \fM_h \ar@{=}[u]
}
\end{equation*}
Since  both $U$ and $W$ are hyperbolic, $\xi$ is rigid by Lemma \ref{lem:rigidity_maps_hyperbolic_stack_curves}. That is,   since $\xi'$ is a deformation of $\xi$, $\xi = \xi'_t $ for every $t \in T$. However, then also  $\nu'_t = \nu$ for every $t \in T$. Therefore, $\nu$ is rigid, indeed.
\end{proof}

\bibliographystyle{skalpha}
\bibliography{includeNice}
 
\end{document}